\documentclass{fic-l}
\usepackage{amssymb}
\newtheorem{theorem}{Theorem}[section]
\newtheorem{corollary}[theorem]{Corollary}
\newtheorem{lemma}[theorem]{Lemma}
\newtheorem{proposition}[theorem]{Proposition}
\newtheorem{example}[theorem]{Example}
\newtheorem{setting}[theorem]{Setting}
\newtheorem{conclusion}[theorem]{Conclusion}

\theoremstyle{definition}
\newtheorem{definition}[theorem]{Definition}

\theoremstyle{remark}
\newtheorem{remark}[theorem]{Remark}

\numberwithin{equation}{section}

\newcommand{\bA}{{\mathbb A }}
\newcommand{\CC}{{\mathbb C }}
\newcommand{\FF}{{\mathbb F }}
\newcommand{\GG}{{\mathbb G }}
\newcommand{\HH}{{\mathbb H }}
\newcommand{\KK}{{\mathbb K }}
\newcommand{\NN}{{\mathbb N }}
\newcommand{\PP}{{\mathbb P }}
\newcommand{\QQ}{{\mathbb Q }}

\newcommand{\ZZ}{{\mathbb Z }}

\newcommand{\XX}{{\mathbb X }}
\newcommand{\bS}{{\mathbb S }}
\newcommand{\cO}{{\mathcal O }}
\newcommand{\cC}{{\mathcal C }}

\newcommand{\cJ}{{\mathcal J }}
\newcommand{\cP}{{\mathcal P }}
\newcommand{\wtF}{{\Phi}}

\newcommand{\wom}{\underline{\mathcal {W}\Omega}{}}
\newcommand{\wo}{\underline{\mathcal {WO}}{}}
\newcommand{\WA}{\underline{{\mathcal {W}}A}}

\newcommand{\half}{[\textstyle{\frac{1}{2}}]}

\newcommand{\wan}[1]{\underline{\mathcal {W}}_{#1}\underline{A}}
\newcommand{\won}[1]{\underline{\mathcal {W}}_{#1}\underline{\mathcal{O}}}
\newcommand{\womn}[1]{\underline{\mathcal {W}}_{#1}\underline{\Omega}}
\newcommand{\q}[1]{\{#1\}}

\newcommand{\tomega}{\widetilde{\Omega}{}}
\newcommand{\OO}{\mathcal {O}}
\newcommand{\oF}{\overline{F}}

\newcommand{\OM}{\Omega}

\newcommand{\Fil}{\mathrm{Fil}}

\newcommand{\FH}{\mathrm{Fil}_\mathrm{Hodge}}
\newcommand{\FC}{\mathrm{Fil}^\mathrm{con}}

\newcommand{\rd}{\mathsf{d}}

\newcommand{\hw}{\mathfrak{h}}

\newcommand{\nil}{\mathcal{N}}
\newcommand{\nar}{\mathfrak{Nilalgs}_{\bA}}

\newcommand{\Spec}{\mathrm{Spec}\:}

\newcommand{\dul}[1]{\underline{\underline{#1}}}

\newcommand{\tr}{{\mathrm{1\hspace{-.25em}l}}}

\begin{document}

\title{The Ordinary Limit for Varieties over $\ZZ[x_1,\ldots,x_r]$}
\author{Jan Stienstra}
\address{Mathematisch Instituut, Universiteit Utrecht\\
Postbus 80.010,  3508 TA Utrecht, the Netherlands}
\email{stien@math.uu.nl}
\begin{abstract}
We investigate for families of smooth projective varieties over a localized
polynomial ring
$\ZZ[x_1,\ldots,x_r][D^{-1}]$ the conjugate filtration on De Rham cohomology
$\otimes\ZZ/N\ZZ$. As $N$ tends to $\infty$
this leads to the concept of the
\emph{ordinary limit}, which seems to be the \emph{non-archimedean analogue} of
the \emph{large complex structure limit}.
\end{abstract}

\maketitle

\section*{Introduction}

Many important families of Calabi-Yau threefolds appear in the following form:
There is a ring $\bA=\ZZ[x_1,\ldots,x_r][D^{-1}]$, with $D$ some polynomial in
the polynomial ring $\ZZ[x_1,\ldots,x_r]$, and
there is a smooth projective morphism
$f:\XX\rightarrow\bS=\Spec \bA$ of relative dimension $3$
such that all De Rham cohomology groups $H_{DR}^m(\XX/\bS)$ and all Hodge
cohomology groups $H^j(\XX,\OM_{\XX/\bS}^i)$ are free $\bA$-modules and
such that
\begin{equation}\label{eq:CY condition}
\OM_{\XX/\bS}^3\simeq \cO_{\XX}\,,\quad H^1(\XX,\cO_\XX)=H^2(\XX,\cO_\XX)=0\,.
\end{equation}

This happens, for instance, for complete intersections in a projective space
$\PP^d_\bA$ given by homogeneous polynomials with coefficients in $\bA$ and
with degrees summing to $d+1$, like $\diamond$ quintic hypersurfaces in $\PP^4$
$\diamond$ intersections of two cubics in $\PP^5$ $\diamond$ intersections of
four quadrics in $\PP^7$.
The conditions of smoothness of $\XX/\bS$ and freeness of cohomology then boil
down to a condition that certain polynomial expressions in the coefficients
should be invertible in $\bA$.

In connection with Mirror Symmetry one is particularly interested in solutions
of the Picard-Fuchs equations of the given family of Calabi-Yau threefolds. The
Picard-Fuchs equations describe how a non-zero global differential $3$-form
varies in the family. For families like those in the first paragraph the
Picard-Fuchs equations are completely defined over $\bA$
and there are many analytical environments in which one may look for solutions.
Traditionally one works in an \emph{archimedean environment} represented by
complex geometry and complex functions and one looks at Calabi-Yau threefolds
\emph{near the large complex structure limit}.
In the present paper we want to work in a \emph{non-archimedean environment}
which is represented by projective systems of groups
$\left\{M_N \right\}_{N\in\NN}$, with $M_N$ a module over $\ZZ/N\ZZ$, indexed
by the positive integers with their divisibility relation, i.e. for every
$N\in\NN$ there is given a $\ZZ/N\ZZ$-module $M_N$ and for every pair of
positive integers $(N,K)$ with $K$ dividing $N$ there is given a homomorphism
$M_N\rightarrow M_K$.

In \cite{S1} it was pointed out that in the crystalline cohomology of
\emph{families of ordinary Calabi-Yau threefolds in the $p$-adic setting}
the structure is analogous to that for Calabi-Yau threefolds near the
large complex structure limit. The analogy comes from the fact that in both the
complex and the $p$-adic situation there is a filtration on De Rham cohomology
opposite to the Hodge filtration, stable under the Gauss-Manin connection and
with associated graded module of Hodge-Tate type.
For complex Calabi-Yau threefolds near the large complex structure limit
one obtains the filtration from the \emph{local monodromy around this limit
point} \cite{D,M}.
For families of ordinary Calabi-Yau threefolds over a base of positive
characteristic one uses the conjugate filtration and the action of
\emph{Frobenius operators}.
The present paper does not build on \cite{S1}, but takes only the suggestion
that by studying ``ordinariness'' and Frobenius operators one may discover
interesting arithmetic-geometrical facts about Calabi-Yau threefolds.
Ordinariness is a well defined condition for varieties in characteristic
$p>0$.
The present paper focusses on the question: \emph{What is the right notion of
``ordinariness'' for varieties over $\ZZ[x_1,\ldots,x_r]$?}
For reasons of time and space we must postpone the discussion of topics which
are specific for Calabi-Yau threefolds, like \emph{canonical coordinates} and
\emph{Yukawa coupling} to another article. Without specific CY3 conditions the
setting becomes:

\

\begin{setting}\label{setting}
There is given a smooth projective morphism
$f:\XX\rightarrow\bS=\Spec \bA$ of relative dimension $\rd$ with ring $\bA$
\'etale over a polynomial ring $\ZZ\half[x_1,\ldots,x_r]$ and
such that all De Rham cohomology groups $H_{DR}^m(\XX/\bS)$ and all Hodge
cohomology groups $H^j(\XX,\OM_{\XX/\bS}^i)$ are free $\bA$-modules.
The $\bA$-module $H^\rd(\XX,\OM^\rd_{\XX/\bS})$ should have rank $1$ and be
generated by an element $\varpi$ from $H^\rd(\XX,\OM^\rd_{\XX/\bS,\log})$ (cf.
(\ref{eq:log forms})). Finally for all $i,j$ the pairing
$ \langle\:,\:\rangle:\;
H^j(\XX,\OM^i_{\XX/\bS})\times H^{\rd -j}(\XX,\OM^{\rd -i}_{\XX/\bS})
\rightarrow \bA$, defined by the rule
$\alpha\cdot\beta=\langle\alpha,\beta\rangle\varpi$,
should be non-degenerate.
\end{setting}

\

Note that this setting behaves well with respect to \'etale base change:
\\
if $A$ is an \'etale algebra over $\bA$, $S=\Spec A$ and $X=\XX\times_{\bS} S$,
then $f:X\rightarrow S$ enjoys the same properties;
namely
$H_{DR}^m(X/S)=H_{DR}^m(\XX/\bS)\otimes_\bA A$ for every $m$ and
$H^j(X,\OM_{X/S}^i)=H^j(\XX,\OM_{\XX/\bS}^i)\otimes_\bA A$ for all $i,j$
and these
are free $A$-modules; moreover, the induced pairing
$
H^j(X,\OM^i_{X/S})\times H^{\rd -j}(X,\OM^{\rd -i}_{X/S})
\rightarrow A$ is non-degenerate.

The reason for requiring that $2$ be invertible in $\bA$ is some technical
$2$-torsion problem in the constructions of \cite{S2}; see loc. cit. 2.7-2.9.
The reason for having $\bA$ \'etale over $\ZZ\half[x_1,\ldots,x_r]$ instead of
just of the form $\ZZ[x_1,\ldots,x_r][D^{-1}]$ is that for a continuation of
this work one needs ``vectors fixed by Frobenius''. That requires solving
certain systems of polynomial equations and thus more involved \'etale
extensions than just Zariski localizations will be needed. Therefore we choose
a formulation so that it is clear that the results remain valid after such an
extension.

For varieties in characteristic
$p>0$ ``ordinariness'' can be described in various equivalent ways. In Section
\ref{section:conjugate} we look at the condition of ordinariness which in
characteristic $p$ requires the \emph{vanishing of all cohomology groups of all
sheaves of exact differential forms}. Its straightforward generalization to our
setting can be formulated as: \textit{for fixed} $N\geq 2$:
\begin{equation}\label{eq:ZN iso intro}
H^j(X,Z_N\OM_{X/S}^i)\rightarrow H^j(X,\OM_{X/S}^i)
\quad\textit{is an isomorphism for all}\quad i,j,
\end{equation}
where
$Z_N\OM^i_{X/S}:=\{\omega\in\OM^i_{X/S}\:|\:d\omega\in N\cdot\OM^{i+1}_{X/S}\}$
is the sheaf of $i$-forms which are closed modulo $N$.
Here we write $X/S$ and not $\XX/\bS$ because the condition will usually only
be satisfied after we restrict the original family to a (Zariski) open subset
$S$ of $\bS$ which depends on $N$.
Theorem \ref{thm:conj split} and Corollary \ref{cor:G-M on limit} state that if
Condition
(\ref{eq:ZN iso intro}) is satisfied the \emph{conjugate filtration} on
$H^m_{DR}(X/S)\otimes \ZZ/N\ZZ$, for every $m$, is indeed opposite to
the Hodge filtration and is stable for the Gauss-Manin connection.
That is, however, not all we want. For one thing, as $N$ varies, conditions and
results should fit into projective systems. For the condition (\ref{eq:ZN iso
intro}) this can be achieved (see Proposition \ref{prop:reduce ordinary})
by replacing it by the (stronger) condition
\begin{equation}\label{eq:Zp iso intro}
\begin{array}{rl}
H^j(X,Z_p\OM_{X/S}^i)\rightarrow H^j(X,\OM_{X/S}^i)
&\textit{is an isomorphism for all }\; i,j\;
\textit{ and}\\
&\textit{all prime numbers }\;p\;\textit{ dividing }\; N.
\end{array}
\end{equation}
More importantly
we also want \emph{the associated graded module to be of Hodge-Tate type}.
The associated graded module is the \emph{Hodge cohomology} $\bmod N$:
$$
\bigoplus_{i+j=m} H^j(X,\OM^i_{X/S})\otimes\ZZ/N\ZZ.
$$
To formulate what ``Hodge-Tate type'' means we need Frobenius operators.
The traditional point of view on Frobenius operators is that they arise in
characteristic $p$, for only one prime number $p$ at a time.
Here we want to work at the level of algebraic geometry over
$\ZZ\half[x_1,\ldots,x_r]$. In Sections
\ref{section:FrobWitt}--\ref{section:conjugate2} we describe the formalism of
\emph{generalized Witt vectors} and the \emph{generalized De Rham-Witt complex}
which yields a Frobenius operator and a Verschiebung operator for every
positive integer. In Section
\ref{section:ZN frob} we relate $Z_N\OM^i_{X/S}$
to Frobenius and Verschiebung operators.
In Theorem \ref{thm:U=conj} we relate the conjugate filtration, the Hodge
filtration and the Hodge decomposition to what we call the
\emph{Hodge-Witt cohomology of $X/S$}. In Section \ref{section:conjugate3}
we investigate when the Hodge-Witt cohomology of $X/S$ is of ``Hodge-Tate
type''.
It turns out that the conditions must be stengthened once more: instead of
taking in (\ref{eq:Zp iso intro}) only primes $p$ dividing $N$ one should take
all primes $p\leq N$. Thus we are led to propose

\begin{definition}\label{def:ordinary up to N}
Let $f:X\rightarrow S=\Spec A$ satisfy the hypotheses in Setting \ref{setting}.
Let $N$ be a positive integer.
We say that \emph{$X/S$ is ordinary up to level $N$} if
\begin{equation}\label{eq:ordinary up to N}
\begin{array}{rl}
H^j(X,Z_p\OM_{X/S}^i)\rightarrow H^j(X,\OM_{X/S}^i)
&\textit{is an isomorphism for all }\; i,j\;
\textit{ and}\\
&\textit{all prime numbers }\;p\leq N.
\end{array}
\end{equation}
\end{definition}

\

In Theorem \ref{thm:Hodge-Tate} we show that if $X/S$ is ordinary up to level
$N$, then the Frobenius operators on Hodge-Witt cohomology induce for every
prime number $p\leq N$ and for all $i,j$ an isomorphism
\begin{equation}\label{eq:HT}
\oF_p: \oF_p^*H^j(X,\OM_{X/S}^i)\stackrel{\simeq}{\rightarrow}
H^j(X,\OM_{X/S}^i)\otimes\ZZ/p\ZZ
\end{equation}
where
$$
\oF_p^*H^j(X,\OM_{X/S}^i)=(A/pA)\otimes_A H^j(X,\OM_{X/S}^i)
$$
with $A/pA$ viewed as an $A$-module via the ring homomorphism
$A\rightarrow A/pA\,,\quad a\mapsto a^p\bmod p$. Isomorphisms like
(\ref{eq:HT}) induced by Frobenius operators are the standard criterion for
being of Hodge-Tate type. In the present situation it is probably better to say
\emph{of Hodge-Tate type up to level $N$}.

\

Setting \ref{setting}, Definition \ref{def:ordinary up to N},
Theorems \ref{thm:conj split}, \ref{thm:U=conj} and \ref{thm:Hodge-Tate}
and the above description of ``Hodge-Tate type up to level $N$'' lead to the
following conclusion.

\begin{conclusion}\label{conclusion}
Let $f:\XX\rightarrow \bS=\Spec \bA$ be as in Setting \ref{setting}.
Let $A$ be an \'etale $\bA$-algebra, $S=\Spec A$ and $X=\XX\times_\bS S$.
Assume that $X/S$ is ordinary up to level $N$, for some positive integer $N$.
Then for every $m$ and for every $n\leq N$ the conjugate filtration
on $H^m_{DR}(X/S)\otimes\ZZ/n\ZZ$ is opposite to the Hodge filtration,
is stable for the Gauss-Manin connection and the associated graded
$\bigoplus_{i+j=m} H^j(X,\OM^i_{X/S})\otimes\ZZ/n\ZZ$ is of
Hodge-Tate type up to level $N$.
\qed
\end{conclusion}

\

This conclusion gives a somewhat weakened description of the actual structure.
For a more precise description Propositions
\ref{prop:PhiN compatibilities}, \ref{prop:wittsurjective1} and Corollary
\ref{cor:wittsurjective i} should also be taken into account. Moreover these
suggest a reformulation in terms of \emph{formal groups}, which might turn out
to be quite attractive. A hint to this formal group structure is given in the
Appendix.

\

The condition that $X/S$ be ordinary up to level $N$ will in general only be
satisfied after restricting the original family $\XX/\bS$ to an open subset $S$
of $\bS$. This subset depends on $N$. As $N$ moves up through $\NN$ the
corresponding open set will shrink. In the limit we have a subset
$S_\infty\subset \bS$.
\begin{itemize}
\item
\emph{We want to call $S_\infty$ the ordinary limit set of $\XX/\bS$}.
\end{itemize}
$S_\infty$ contains $\Spec(\bA\otimes\QQ)$, because there the conditions for
ordinariness up to any level are trivially satisfied. But over
$\Spec(\bA\otimes\QQ)$ the De Rham cohomology groups
modulo $N$ are just $0$. So there there is no interesting conclusion.
However $S_\infty$ should be much larger than $\Spec(\bA\otimes\QQ)$
and should contain an affine set $\Spec(A_\infty)$ on which no or just a few
primes are invertible; see the example discussed below and also Remark
\ref{remark:ordinary open}.
The geometry of the family $\XX$ restricted to $\Spec(A_\infty)$ is not yet
clear (to me), but examples like the one below indicate that
$\Spec(A_\infty)$ contains a punctured formal neighborhood of the
``large complex structure limit point'' and that the
aforementioned formal groups and their interaction via the Gauss-Manin
connection extend over this limit point.
I like to view this as analogous to the traditional complex situation
where the limit fibre is a singular Calabi-Yau threefold and where one has a
limit mixed Hodge structure.

\

To illustrate some of the above issues let us look at the pencil of elliptic
curves
\begin{equation}\label{eq:hesse}
x(X^3+Y^3+Z^3)+XYZ=0,
\end{equation}
with the nine base points blown up.
There are singular fibres for $x=0$ and $x^3=\frac{1}{27}$.
Thus we work over the ring
$\bA=\ZZ[x,\frac{1}{2x(27 x^3-1})]$.

The formal group law (in an appropriate coordinatization) for these
elliptic curves is
\begin{equation}\label{eq:group law}
G(t_1,t_2)=\ell^{-1}(\ell(t_1)+\ell(t_2)).
\end{equation}
with
$$
\ell(t)\:=\:\sum_{m\geq 1}
\frac{1}{m}\,a_m(x)t^m\;\in (\bA\otimes\QQ)[[t]]
$$
and
\begin{equation}\label{eq:hessecoeff}
a_m(x)=\sum_{j\geq 0} \frac{(3j)!}{j!^3}
\left(\begin{array}{c}m-1\\ 3j\end{array}\right) x^{3j}\,\in\ZZ[x].
\end{equation}
For $\lambda$ in some field of characteristic $p>2$ it is well known
that the elliptic curve
$\lambda(X^3+Y^3+Z^3)+XYZ=0$ is ordinary if and only if
its \emph{Hasse-Witt invariant} $a_p(\lambda)$ is not zero; here we
use one of many equivalent characterizations of ordinariness for elliptic
curves in characteristic $p$.
So the ordinary limit set $S_\infty$ in this example is a union of
affine sets $\Spec A_\Pi$ where $\Pi$ runs over all (finite as well
as infinite) subsets of the set $\cP$ of all prime numbers and
$$
A_\Pi := \bA[a_p(x)^{-1},\,r^{-1}\:|\:p\in\Pi,\:r\in\cP\setminus\Pi ].
$$
Thus, for $\Pi=\emptyset$, the empty set, $A_\emptyset=\bA\otimes\QQ$.
On the other extreme, $\Pi=\cP$ gives the ring we have in mind for the above
mentioned $A_\infty$; so,
$$
A_\infty=A_\cP=\bA[a_p(x)^{-1}\:|\:p\;\textrm{prime}].
$$
Note that $a_p(0)=1$ for all $p$ and that, hence, $A_\infty$ embeds into
the Laurent series ring $\ZZ[\frac{1}{2}][[x]] [x^{-1}]$. Thus we see that
$\Spec A_\infty$ contains a punctured formal neighborhood of $x=0$.
The fibre over $x=0$ is singular, consisting of three lines, and the monodromy
around this fibre is maximally unipotent. In the physicists' language $x=0$ is
the large complex structure limit point.
Note that $a_m(0)=1$ for all $m\geq 1$ and that, hence, the formal group
law (\ref{eq:group law}) extends well over $x=0$; in fact at $x=0$ it is the
standard multiplicative group law $t_1+t_2-t_1t_2$.

For families of K3-surfaces, like  $x(W^4+X^4+Y^4+Z^4)+WXYZ=0$, the story is
the same (with $3$ replaced by $4$ in (\ref{eq:hessecoeff})). In particular the
fibre at $x=\lambda$ is ordinary if its Hasse-Witt invariant is not zero. The
relevant formal group is the formal Brauer group $H^2(\XX,\hat{\GG}_{m,\XX})$;
see the appendix.

For families of Calabi-Yau threefolds invertibility of the Hasse-Witt invariant
associated with the Artin-Mazur formal group $H^3(\XX,\hat{\GG}_{m,\XX})$ (see
the appendix)
is necessary, but
may be not sufficient to ensure ordinariness. Yet I do expect that also in this
CY3 case the ordinary limit set contains a punctured formal neighborhood of
``the large complex structure limit point''.

\section{The conjugate filtration and ordinariness: Act 1}
\label{section:conjugate}
\subsection{Definition and basic properties of the conjugate filtration.}
For a smooth projective morphism $X\rightarrow S$ of schemes the De Rham
cohomology
$H_{DR}^m(X/S)$ is by definition the hypercohomology
$\HH^m(X,\OM_{X/S}^\bullet)$ of the De Rham complex
\footnote{all sheaves in this paper are taken with respect to the Zariski
topology}
$$
\OM_{X/S}^\bullet\,:\; \ldots \rightarrow \cO_X
\stackrel{d}{\rightarrow}
\OM_{X/S}^1\stackrel{d}{\rightarrow}\OM_{X/S}^2
\stackrel{d}{\rightarrow}\ldots\stackrel{d}{\rightarrow}
\OM_{X/S}^{i-1}\stackrel{d}{\rightarrow}
\OM_{X/S}^i\stackrel{d}{\rightarrow}
\OM_{X/S}^{i+1}\stackrel{d}{\rightarrow}\ldots
$$
Every complex $\cC^\bullet$ carries two natural filtrations by sub-complexes
$\cC^{\bullet\geq i}$ and $t_{\leq i}\cC^\bullet$:
\begin{eqnarray}
\label{eq:stupid filtration}
&&
\left(\cC^{\bullet\geq i}\right)^j=0\quad\textrm{if }\; j<i\,,
\qquad \left(\cC^{\bullet\geq i}\right)^j=\cC^j\quad\textrm{if }\;j\geq i.\\
\nonumber
&&
\left(t_{\leq i}\cC^\bullet\right)^j=0\quad\textrm{if }\; j>i\,,
\qquad \left(t_{\leq i}\cC^\bullet\right)^j=\cC^j\quad\textrm{if }\; j<i\,,\\
\label{eq:canonical filtration}
&& \left(t_{\leq i}\cC^\bullet\right)^i=Z\cC^i:=\ker(d:\cC^i\rightarrow
\cC^{i+1}).
\end{eqnarray}
When applied to the De Rham complex $\OM_{X/S}^\bullet$ these induce on
$H_{DR}^m(X/S)$: \\
the \emph{Hodge filtration}
$$
\FH^i H_{DR}^m(X/S)=
\mathrm{image}(\HH^m(X,\OM_{X/S}^{\bullet\geq i})\rightarrow
\HH^m(X,\OM_{X/S}^\bullet))
$$
and the \emph{conjugate filtration}
$$
\FC_i H_{DR}^m(X/S)=
\mathrm{image}(\HH^m(X,t_{\leq i}\OM_{X/S}^\bullet)\rightarrow
\HH^m(X,\OM_{X/S}^\bullet)).
$$
The Hodge filtration is decreasing and the conjugate filtration is increasing
and moreover
\begin{eqnarray*}
&&\FH^0 H_{DR}^m(X/S)\;=\;\FC_m H_{DR}^m(X/S)\;=\;H_{DR}^m(X/S)\,,
\\
&&\FH^{m+1} H_{DR}^m(X/S)\;=\;\FC_{-1} H_{DR}^m(X/S)\;=\; 0\,.
\end{eqnarray*}
Let
$B\OM_{X/S}^i:=
\mathrm{image}(d:\OM_{X/S}^{i-1}\rightarrow\OM_{X/S}^i)$
resp.
$Z\OM_{X/S}^i:=
\mathrm{ker}(d:\OM_{X/S}^i\rightarrow\OM_{X/S}^{i+1})$ denote
the sheaves of exact resp. closed $i$-forms.
There are short exact sequences of complexes
\begin{eqnarray}\label{eq:conj-hodge sequence}
&&0\rightarrow t_{\leq i}\OM_{X/S}^\bullet\oplus\OM_{X/S}^{\bullet\geq i+1}
\rightarrow \OM_{X/S}^\bullet\rightarrow B\OM_{X/S}^{i+1}[-i]\rightarrow 0
\\
&&0\rightarrow Z\OM_{X/S}^i[-i]\rightarrow t_{\leq
i}\OM_{X/S}^\bullet\oplus\OM_{X/S}^{\bullet\geq i}
\rightarrow \OM_{X/S}^\bullet\rightarrow 0
\\
&&0\rightarrow Z\OM_{X/S}^i\rightarrow \OM_{X/S}^i\rightarrow
B\OM_{X/S}^{i+1}\rightarrow 0
\end{eqnarray}
where $B\OM_{X/S}^{i+1}[-i]$ ( resp. $Z\OM_{X/S}^i[-i]$) is the
complex with
$B\OM_{X/S}^{i+1}$ ( resp. $Z\OM_{X/S}^i$) sitting in degree $i$ and all other
terms equal to $0$. The corresponding exact sequences of hypercohomology groups
show: \emph{if the condition}
\begin{equation}
\label{eq:ordinarycondition}
H^{m-i}(X,B\OM_{X/S}^{i+1})=H^{m-i-1}(X,B\OM_{X/S}^{i+1})=
H^{m-i}(X,B\OM_{X/S}^i)=0
\end{equation}
\emph{is satisfied, then}
\begin{eqnarray}
\label{eq:opposite1}
H_{DR}^m(X/S)&=&\FC_i H_{DR}^m(X/S)\oplus\FH^{i+1} H_{DR}^m(X/S)
\\[1ex]
\label{eq:opposite2}
H^{m-i}(X,\OM_{X/S}^i)&=&\FC_i H_{DR}^m(X/S)\cap\FH^i H_{DR}^m(X/S).
\end{eqnarray}
So if (\ref{eq:ordinarycondition}) holds for
$i=0,\ldots,m$, then the Hodge filtration and the conjugate filtration on
$H_{DR}^m(X/S)$ are opposite and one has the \emph{Hodge decomposition}:
\begin{equation}\label{eq:hodgedecomposition}
H_{DR}^m(X/S)=\bigoplus_{i=0}^mH^{m-i}(X,\OM_{X/S}^i).
\end{equation}

\begin{remark}
This form of conjugate filtration plays no role in complex geometry, since for
varieties over $\CC$ with the complex topology the \emph{Poincar\'e lemma}
says that locally all closed differential forms are exact, and thus implies
that the filtration $\{\FC_i\}_{}$ jumps from null at $i=-1$ to all
at $i=0$. On the other hand, according to Hodge theory one obtains for complex
K\"ahler manifolds a
filtration opposite to the Hodge filtration by simply taking the complex
conjugate of the Hodge filtration.
\end{remark}

Condition (\ref{eq:ordinarycondition}) is a familiar condition for smooth
projective varieties in positive characteristic \cite{I2,IR}:
\textit{When $f:X\rightarrow S$ is a proper smooth morphism between schemes
of characteristic $p>0$, one says that $X$ is \emph{ordinary over} $S$ if
for all $i,j$}
$$
R^j f_* B\OM^i_{X/S}=0.
$$
In case $S=\Spec k$ with $k$ a perfect field  of characteristic $p>0$
this is one of many equivalent ways to
define ordinariness (cf. \cite{IR} thm. IV 4.13). This one involves no concepts
specific for characteristic $p$ and it is tempting to use it also in other
contexts
as a quick way for splitting the Hodge filtration, with the conjugate
filtration as its opposite.
We want to work in the geometric setting of a smooth projective morphism
$f:X\rightarrow S=\Spec A$ in which the ring $A$ is \'etale over a polynomial
ring $\ZZ[x_1,\ldots,x_r]$. We have found that the following makes a good
bridge between the various prime characteristics:
\begin{itemize}
\item
\emph{Work systematically with projective systems of groups
$\{M_N\}_{N\in\NN}$, with $M_N$ a module over $\ZZ/N\ZZ$, indexed by the set
$\NN$ which is ordered by divisibility.}
\end{itemize}
Throughout this paper we use for an integer $N$ and an abelian group $G$ the
\\
\textbf{notation:}
\begin{equation}\label{eq:mod N}
\begin{array}{lll}
G\q{N}&:=& G/NG\quad \textrm{ in case of additive notation} \\
&:=& G/G^N \quad \textrm{ in case of multiplicative notation}.
\end{array}
\end{equation}

Thus instead of $\OM_{X/S}^\bullet$ we take the complex
$\OM_{X/S}^\bullet\q{N}$.
If the De Rham cohomology group $H_{DR}^m(X/S)$ and the Hodge cohomology group
$H^{m-i}(X,\OM_{X/S}^i)$ are free $A$-modules, then
\begin{equation}\label{eq:coho mod N}
\begin{array}{rcl}
\HH^m(X,\OM_{X/S}^\bullet\q{N})&=&H_{DR}^m(X/S)\q{N}\,,\\
H^{m-i}(X,\OM_{X/S}^i\q{N})&=&H^{m-i}(X,\OM_{X/S}^i)\q{N}\,.
\end{array}
\end{equation}
Note that
$$
\OM_{X/S}^\bullet\q{N}=
\OM^\bullet_{X\times\Spec(\ZZ/N\ZZ)\,/\,S\times\Spec(\ZZ/N\ZZ)}\,.
$$
The above arguments can be applied to
$X\times\Spec(\ZZ/N\ZZ)\stackrel{f}{\rightarrow}S\times\Spec(\ZZ/N\ZZ)$.
Before stating our conclusion
we make a closer analysis of the condition
\begin{equation}\label{eq:exact vanish}
H^j(X,B(\OM_{X/S}^i\q{N}))=0
\qquad\textit{for all }\;i,j.
\end{equation}
The cohomology sequences associated with the short exact sequences of complexes
$$
0\rightarrow Z(\OM_{X/S}^i\q{N})\rightarrow \OM_{X/S}^i\q{N}\rightarrow
B(\OM_{X/S}^{i+1}\q{N})\rightarrow 0
$$
show that the condition (\ref{eq:exact vanish}) is equivalent with:
\begin{equation}\label{eq:closed iso}
H^j(X,Z(\OM_{X/S}^i\q{N}))\rightarrow H^j(X,\OM_{X/S}^i\q{N})
\quad\textit{is an isomorphism for all}\quad i,j.
\end{equation}
There is yet another equivalent form of this condition. For that we consider
\begin{equation}\label{eq:ZN}
Z_N\OM^i_{X/S}:=\{\omega\in\OM^i_{X/S}\:|\:d\omega\in N\cdot\OM^{i+1}_{X/S}\}.
\end{equation}
It fits into the commutative diagram with exact rows
$$
\begin{array}{ccccccccc}
0&\rightarrow&\OM^i_{X/S}&\stackrel{\cdot N}{\longrightarrow}&
Z_N\OM^i_{X/S}&\longrightarrow&Z(\OM^i_{X/S}\q{N})&\rightarrow&0\\
&&\parallel&&\downarrow&&\downarrow&&\\
0&\rightarrow&\OM^i_{X/S}&\stackrel{\cdot N}{\longrightarrow}&
\OM^i_{X/S}&\longrightarrow&\OM^i_{X/S}\q{N}&\rightarrow&0
\end{array}
$$
Under the assumptions on $X/S$ made in Theorem \ref{thm:conj split}
multiplication by $N$ on $H^j(X,\OM^i_{X/S})$ is injective. Therefore the
ladder of cohomology groups for the above diagram splits into diagrams with
exact rows
$$
\begin{array}{ccccccccc}
0&\!\rightarrow\!&H^j(X,\OM^i_{X/S})&\!\stackrel{\cdot N}{\rightarrow}\!&
H^j(X,Z_N\OM^i_{X/S})&\!\rightarrow\!&H^j(X,Z(\OM^i_{X/S}\q{N}))
&\!\rightarrow\!&0\\
&&\parallel&&\downarrow&&\downarrow&&\\
0&\!\rightarrow\!&H^j(X,\OM^i_{X/S})&\!\stackrel{\cdot N}{\rightarrow}\!&
H^j(X,\OM^i_{X/S})&\!\rightarrow\!&H^j(X,\OM^i_{X/S}\q{N})&\!\rightarrow\!&0
\end{array}
$$
This shows that condition (\ref{eq:closed iso}), and hence also
condition (\ref{eq:exact vanish}), is equivalent with condition
(\ref{eq:ZN iso}) below.

\

Our conclusion from the above discussion is:

\begin{theorem}\label{thm:conj split}
Let $f:X\rightarrow S=\Spec A$ be a smooth projective morphism in which the
ring $A$ is \'etale over a polynomial ring $\ZZ[x_1,\ldots,x_r]$.
Assume that $H_{DR}^m(X/S)$ and
$H^{m-i}(X,\OM_{X/S}^i)$ are free $A$-modules for all $m,i$. Fix $N\in\NN$ and
assume
\begin{equation}\label{eq:ZN iso}
H^j(X,Z_N\OM_{X/S}^i)\rightarrow H^j(X,\OM_{X/S}^i)
\quad\textit{is an isomorphism for all}\quad i,j.
\end{equation}
Then one has for all $m,i$:
\begin{eqnarray*}
H_{DR}^m(X/S)\q{N}&=&\FC_i \HH^m(X,\OM_{X/S}^\bullet\q{N})\oplus\FH^{i+1}
H_{DR}^m(X/S)\q{N}
\\
H^{m-i}(X,\OM_{X/S}^i)\q{N}&=&\FC_i \HH^m(X,\OM_{X/S}^\bullet\q{N})\cap\FH^i
H_{DR}^m(X/S)\q{N}
\\
H_{DR}^m(X/S)\q{N}&=&\bigoplus_{i=0}^mH^{m-i}(X,\OM_{X/S}^i)\q{N}.
\end{eqnarray*}
\qed
\end{theorem}

\

\begin{remark}
If $K$ is a divisor of $N$ there are canonical maps, induced by reducing $\bmod
K$ on the De Rham complex $\OM^\bullet_{X/S}\q{N}$,
\begin{eqnarray*}
H_{DR}^m(X/S)\q{N}&\rightarrow&H_{DR}^m(X/S)\q{K}\\
\FC_i \HH^m(X,\OM_{X/S}^\bullet\q{N})&\rightarrow&
\FC_i \HH^m(X,\OM_{X/S}^\bullet\q{K})\\
\FH^{i+1} H_{DR}^m(X/S)\q{N}&\rightarrow&
\FH^{i+1} H_{DR}^m(X/S)\q{K}
\end{eqnarray*}
and thus the conclusion part of Theorem \ref{thm:conj split} fits well into a
projective system indexed by the positive integers ordered by divisibility.

It is however not clear that if condition (\ref{eq:ZN iso}) holds for $N$, it
also holds with $K$ in place of $N$. So the condition part of Theorem
\ref{thm:conj split} does not fit well into a projective system.
Proposition \ref{prop:reduce ordinary} will show that this can be remedied by
replacing (\ref{eq:ZN iso}) by the condition
\begin{equation}\label{eq:Zp iso}
\begin{array}{rl}
H^j(X,Z_p\OM_{X/S}^i)\rightarrow H^j(X,\OM_{X/S}^i)
&\textit{is an isomorphism for all }\; i,j\;
\textit{ and}\\
&\textit{all prime numbers }\;p\;\textit{ dividing }\; N.
\end{array}
\end{equation}
\end{remark}

\begin{remark}\label{remark:ordinary open}
Theorem \ref{thm:conj split} is to be used as follows. Let us start with a
smooth projective morphism
$f:\XX\rightarrow \bS=\Spec \bA$ in which the ring $\bA$ is \'etale over a
polynomial ring $\ZZ[x_1,\ldots,x_r]$ and assume that $H_{DR}^m(\XX/\bS)$ and
$H^{m-i}(\XX,\OM_{\XX/\bS}^i)$ are free $\bA$-modules for all $m,i$.
Then, for a fixed $N$, condition (\ref{eq:ZN iso}) will in general not be
satisfied for the whole $\XX/\bS$. Instead there will be a non-empty
open (for the Zariski topology) $S\subset \bS$ such that
(\ref{eq:ZN iso}) does hold for $X=\XX\times_{\bS} S$. As such this is a
trivial statement since over the open set on which $N$ is invertible
condition (\ref{eq:ZN iso}) is trivially satisfied; but over this set the
conclusion of the theorem also trivializes: all groups involved are zero!
The point is however that the open set on which (\ref{eq:ZN iso}) holds is much
bigger: for $N=p$ a prime number and universal families
$\XX\rightarrow\bS$ of complete intersections in projective space Illusie
\cite{I2} shows that (\ref{eq:exact vanish}), and hence also
(\ref{eq:ZN iso}), holds on a non-empty open part of the characteristic $p$
locus in $\bS$. The complement of this open subset in the characteristic $p$
locus is the zero set of some ideal $\hw_p\subset\bA$ with
$p\bA\subsetneq \hw_p$. Thus condition (\ref{eq:ZN iso}) for $N=p$ will be
satisfied on the complement of the zero locus of the ideal $\hw_p$ in $\bS$.
Proposition \ref{prop:reduce ordinary} will show that for general $N$ the open
set on which
(\ref{eq:ZN iso}) holds is at least as large as the intersection of the open
sets for the prime divisors of $N$ i.e. the complement of the zero set of the
ideal $\prod_{p|N}\hw_p$.
\end{remark}

\subsection{Conjugate filtration and Gauss-Manin connection.}
Let us see how the conjugate filtration behaves with respect to the
Gauss-Manin connection. First recall how Katz and Oda \cite{K} constructed the
Gauss-Manin connection from the filtration on $\OM_{X/\ZZ}^\bullet$ formed by
the powers of the graded ideal \begin{equation}\label{eq:JXS}
\cJ^\bullet_{X/S}:= d(f^*\cO_S)\cdot\OM_{X/\ZZ}^\bullet,
\end{equation}
i.e.
$$
\cJ^0_{X/S}:=\cO_X\,,\qquad
\cJ^i_{X/S}:= d(f^*\cO_S)\cdot\OM_{X/\ZZ}^{i-1}\quad\textrm{for }i\geq1.
$$
First one has to notice that, since $X$ is smooth over $S$,
$$
\OM^\bullet_{X/\ZZ}/\cJ^\bullet_{X/S}=\OM^\bullet_{X/S}\,,\qquad
\cJ^\bullet_{X/S}/(\cJ^\bullet_{X/S})^2=\OM^1_{S/\ZZ}\otimes
\OM^{\bullet-1}_{X/S}.
$$
Then one finds \emph{the Gauss-Manin connection}
\begin{equation}\label{eq:Gauss-Manin}
\nabla: \HH^m(X,\OM_{X/S}^\bullet)\rightarrow
\OM^1_{S/\ZZ}\otimes\HH^m(X,\OM_{X/S}^\bullet).
\end{equation}
as the connecting map
$
\HH^m(X,\OM_{X/S}^\bullet)\rightarrow
\HH^{m+1}(X,\OM^1_{S/\ZZ}\otimes\OM^{\bullet-1}_{X/S})
$
in the long exact cohomology sequence of
$$
0\rightarrow \OM^1_{S/\ZZ}\otimes\OM^{\bullet-1}_{X/S}
\rightarrow \OM_X^\bullet/(\cJ^\bullet_{X/S})^2\rightarrow
\OM_{X/S}^\bullet\rightarrow 0.
$$

This construction immediately shows on the one hand that the Hodge filtration
need not be stable under the Gauss-Manin connection, but instead satisfies
\emph{Griffiths transversality}:
\begin{equation}\label{eq: Griffiths}
\nabla\FH^i H^m_{DR}(X/S)\;\subset\;\OM^1_{S/\ZZ}\otimes\FH^{i-1}
H^m_{DR}(X/S)\,,
\end{equation}
and on the other hand that the conjugate filtration is
stable under the Gauss-Manin connection:
\begin{equation}\label{eq: con stable}
\nabla\FC_i H^m_{DR}(X/S)\;\subset\;\OM^1_{S/\ZZ}\otimes\FC_i H^m_{DR}(X/S).
\end{equation}

All this remains valid when we go modulo some positive integer $N$, i.e. when
we pull $X\rightarrow S\rightarrow \Spec \ZZ$ back to
$X\times\Spec \ZZ/N\ZZ\rightarrow S\times\Spec \ZZ/N\ZZ\rightarrow \Spec
\ZZ/N\ZZ$. As a result we find:

\begin{corollary}\label{cor:G-M on limit}
In the situation of Theorem \ref{thm:conj split} the conjugate filtration on
$H_{DR}^m(X/S)\q{N}$ provides a filtration which is opposite to the Hodge
filtration and which is stable under the Gauss-Manin connection.
\qed
\end{corollary}

\

\section{Generalized Witt vectors}
\label{section:FrobWitt}

For a ring\footnote{all rings in this paper are commutative with $1$, unless
explicitly stated otherwise} $A$ one defines its additive group of
\emph{generalized Witt vectors}\footnote{We add the adjective
\emph{generalized} to emphasize that these are not the more common Witt vectors
associated with a prime number $p$; the latter is the \emph{$p$-typical part}
of $\WA$;  cf. \cite{C,H,B}.}
$\WA$ to be the multiplicative group of one-variable formal power series with
constant term $1$ and coefficients in $A$:
\begin{equation}\label{eq:WA}
\WA :=1+tA[[t]].
\end{equation}
The group $\WA$ naturally comes with the decreasing filtration:
\begin{equation}\label{eq:wittfiltA}
\Fil_n \WA := 1+t^{n+1}A[[t]]\,,\qquad n\geq 0.
\end{equation}
If $a$ is an element of $A$ we write $\dul{a}$ for the
power series of $(1-at)^{-1}$ viewed as an element of $\WA$.  Every element of
$1+t^{n+1}A[[t]]$ can be written uniquely as a $t$-adically converging product
 $\prod_{i\geq n+1}(1-a_it^i)^{-1}$  with all
$a_i\in A$.
For $n\in \NN$ the substitution  $t \mapsto t^n$ induces an (additive)
endomorphism $V_n$ of $\WA$. Thus the elements of $\Fil_n \WA$ can
be written uniquely as a sum
\begin{equation}\label{eq:wvcA}
  \sum_{i\geq n+1}V_i \dul{a_i}
\end{equation}
which converges with respect to the filtration (\ref{eq:wittfiltA}).
The operators $V_n$ are called \emph{Verschiebung operators}.

One can construct a continuous product on $\WA$ so that $\WA$ becomes
a topological ring, and each $\Fil_n \WA$
is an ideal. One can also construct
continuous endomorphisms $F_n$ for $n\in \NN$ (see \cite{C,H}). These are
called \emph{Frobenius operators}. For computations in $\WA$ the following
relations plus continuity with respect to the filtration (\ref{eq:wittfiltA})
suffice.
For $k,m, n\in\NN$, elements $\alpha, \beta\in\WA$ and
$a, b\in A$:
\begin{equation}\label{eq:FVW}
\begin{array}{l}
   \dul{a}\cdot\dul{b} = \dul{ab} \, ,\quad
   F_n\dul{a} = \dul{a}^n\, ,\\[.5ex]
   F_n(\alpha \beta)=(F_n\alpha)(F_n\beta) \, , \quad
   V_n(\alpha(F_n\beta))=(V_n\alpha)\beta\, ,  \\
    F_mV_m=m \, , \quad
   F_mF_n=F_{mn}\, , \quad
   V_mV_n=V_{mn}\, ,\\
   V_kF_m=F_mV_k  \qquad\textrm{if} \quad (k,m)=1.
\end{array}
\end{equation}
The third line of (\ref{eq:FVW}) shows
$\;
F_1=V_1=\textit{ identity operator}.
$
The first line of (\ref{eq:FVW}) shows that the map
\begin{equation}\label{eq:teichmueller}
 A\rightarrow\WA\,,\qquad a\mapsto\dul{a}
\end{equation}
is a \emph{homomorphism of multiplicative monoids}. One calls it the
\emph{Teichm\"uller lifting}.

As an exercise in computing with the relations one can check
\begin{equation}\label{eq:VFfilt}
V_m\Fil_n \WA\subset\Fil_{mn+m-1} \WA\,,\qquad
F_m\Fil_{mn} \WA\subset\Fil_n \WA.
\end{equation}
The ring of \emph{generalized Witt vectors of length $n$} is, by definition,
$$
\wan{n} := \WA / \Fil_n \WA  .
$$
Note
$$
\wan{1}\simeq A\,,\qquad  \dul{a}\bmod\Fil_1 \WA\;\leftrightarrow a.
$$
Since obviously $\Fil_n \WA\subset\Fil_m \WA$ if $n\geq m$ there are standard
\emph{truncation maps}
\begin{equation}\label{eq:trunc}
\tr^n_m: \wan{n}\rightarrow \wan{m}\,.
\end{equation}
When source and target of a truncation map are clear from the context we simply
write $\tr$ instead of $\tr^n_m$.

\

\section{The generalized De Rham-Witt complex}
\label{section:DRW}
The constructions in Section \ref{section:FrobWitt} are functorial in $A$ i.e.
a ring homomorphism $A\rightarrow A'$ induces a ring homomorphism
$\WA\rightarrow\WA'$ compatible with all the truncations, Frobenius and
Verschiebung operators. So one can sheafify the constructions
and thus obtain on every scheme $X$ the \emph{sheaves $\wo_X$ and $\won{n}_X$
of generalized Witt vectors (of length $n$)}, together with
all the truncations, Frobenius and Verschiebung operators.

In \cite{S2} we constructed for every scheme $X$ on which $2$ is invertible
the \emph{generalized De Rham-Witt complex}\footnote{Here the adjective
\emph{generalized} emphasizes that this is not the De Rham-Witt complex in
characteristic $p$ \cite{I}; the latter can be recovered as the
\emph{$p$-typical part} of $\wom^\bullet_X$.}. This is a sheaf
$\wom^\bullet_X$ of anti-commutative differential graded algebras on $X$  with
the following properties (\ref{eq:drw0})-(\ref{eq:FVrel}).

Let $\wom_X^i$ denote the degree $i$ component of
$\wom_X^\bullet$. Then $\wom_X^i = 0$ for $i<0$,
\begin{equation}\label{eq:drw0}
\wom_X^0=\wo_X.
\end{equation}
For every integer $m\geq 1$  and every $i\geq 0$ there are homomorphisms
of additive groups
\begin{equation}\label{eq:FV}
F_m,V_m : \wom_X^i \rightarrow \wom_X^i
\end{equation}
The operators $F_m$ are called \emph{Frobenius operators} and the operators
$V_m$ are called \emph{Verschiebung operators}. On the sheaf of generalized
Witt vectors $\wo_X$  they coincide with the
earlier defined Frobenius and Verschiebung operators.
The following relations hold for all $m, n$, for all sections $\alpha$, $\beta$
of $\wom_X^\bullet$ and all sections $a$ of $\OO_X$:
\begin{equation}\label{eq:FVrel}
\begin{array}{l}
F_nd\dul{a} = \dul{a}^{n-1}d\dul{a}\, ,\\
 F_mV_m=m \,,\quad\quad
   F_mF_n=F_{mn}\,,\quad\quad
   V_mV_n=V_{mn},  \\
 V_md=mdV_m\,,\quad
   dF_m=mF_md\,,\qquad
   F_mdV_m=d ,\\
 F_m(\alpha\beta)=(F_m\alpha)(F_m\beta) \,,\qquad
V_m(\alpha(F_m\beta))=(V_m\alpha)\beta,\\
 V_nF_m=F_mV_n \qquad\textrm{if} \quad (n,m)=1;
\end{array}
\end{equation}
here  $d:\wom_X^i \rightarrow \wom_X^{i+1}$  is the differential of the
differential graded algebra $\wom_X^\bullet$.

\begin{remark}
The relation
$
dF_m=mF_md
$
means on the one hand that $F_m$ does not commute with $d$, but shows on the
other hand that we get an operator $\FF_m$ on $\wom^\bullet_X$ which does
commute with $d$ by defining
\begin{equation}\label{eq:FDRW}
\FF_m=m^i F_m\qquad\textrm{on }\;\wom^i_X.
\end{equation}
For De Rham-Witt complex in characteristic $p$ this is a standard construction,
which is of great importance for the so-called \emph{slope spectral sequence};
see \cite{I,IR}.
\end{remark}

\

The filtration $\{\Fil_n \wo_X\}_{n\geq 0}$ can be extended to a
decreasing filtration on $\wom^\bullet_X$ by defining:
\begin{equation}\label{eq:filtdrw}
\Fil_n\wom^\bullet_X:=\textit{ ideal generated by }\;\Fil_n \wo_X\;
\textit{ and }\;d\Fil_n \wo_X.
\end{equation}
Each $\Fil_n\wom^\bullet_X$ is then a differential graded ideal in
$\wom^\bullet_X$.
The quotient
\begin{equation}\label{eq:drwlevel}
\womn{n}^\bullet_X:=\wom^\bullet_X/\Fil_n\wom^\bullet_X
\end{equation}
is called the \emph{generalized De Rham-Witt complex of level $n$}.
This is a differential graded algebra; $\womn{n}^i_X=0$ if $i<0$ and
$$
\womn{n}^0_X=\won{n}_X.
$$

In particular $\womn{1}^\bullet_X$ is a differential graded algebra with degree
$0$ term equal to $\cO_X$. Let $\Omega_X^\bullet=\Omega_{X/\ZZ}^\bullet$ be the
absolute De Rham complex on $X$ i.e. the complex of differential forms relative
to $\ZZ$. The universal property of $\OM^\bullet_X$ yields a homomorphism of
complexes
$$
\OM^\bullet_X\rightarrow\womn{1}^\bullet_X
$$
which happens to be surjective by \cite{S2} \S 2.16.
On the other hand it has been shown in \cite{S2} \S 2.20 that there is
a homomorphism of sheaves of differential graded algebras
$$
\wom_X^\bullet \rightarrow \tomega_X^\bullet;
$$
where
$ \tomega_X^i := \Omega_X^i /(i!\mbox{-torsion in } \Omega_X^i)$.
The composite of these two maps is the obvious map
$\OM^\bullet_X\rightarrow \tomega_X^\bullet$. \emph{Very mild conditions on $X$
will remove the $i!$-torsion and thus assure that}
\begin{equation}\label{eq:dr=drw1}
\womn{1}^\bullet_X\,=\,\OM^\bullet_X.
\end{equation}

\

Since obviously $\Fil_n \wom^\bullet_X\subset\Fil_m \wom^\bullet_X$ if $n\geq
m$ there are standard \emph{truncation maps}
\begin{equation}\label{eq:trunc drw}
\tr^n_m: \womn{n}_X\rightarrow \womn{m}_X\,.
\end{equation}
When source and target of a truncation map are clear from the context we simply
write $\tr$ instead of $\tr^n_m$.

\

\begin{lemma}
For all $m$ and $n$ one has
\begin{equation}\label{eq:VFfilt2}
V_m\Fil_n \wom^\bullet_X\subset\Fil_{mn+m-1} \wom^\bullet_X\,,\qquad
F_m\Fil_{mn} \wom^\bullet_X\subset\Fil_n \wom^\bullet_X.
\end{equation}
Consequently $V_m$ and $F_m$ induce maps
\begin{equation}\label{eq:VFfilt3}
V_m:\womn{n}^\bullet_X\rightarrow\womn{nm+m-1}^\bullet_X\,,\qquad
F_m:\womn{nm}^\bullet_X\rightarrow\womn{n}^\bullet_X.
\end{equation}
\end{lemma}
\begin{proof}
Most of this follows easily from the relations (\ref{eq:FVrel}) and
(\ref{eq:VFfilt}). Only for
$F_md\Fil_{mn} \wo_X\subset\Fil_n \wom^1_X$ one needs the following slightly
tricky calculation: take $i>mn$; let $r$ be the greatest common divisor of $i$
and $m$; set $i'=i/r$ and $m'=m/r$; take integers $j,k$ such that
$ji'+km'=1$; then
\begin{equation}\label{eq:trick}
\begin{array}{rcl}
F_mdV_i\dul{a}&=&F_{m'}dV_{i'}\dul{a}\;=\;jF_{m'}V_{i'}d\dul{a}+
kdF_{m'}V_{i'}\dul{a}\\[.5ex]
&=&jV_{i'}(\dul{a}^{m'-1}d\dul{a})+kdV_{i'}(\dul{a}^{m'})
\\[.5ex]
&=&jV_{i'}(\dul{a}^{m'-1})d(V_{i'}\dul{a})+kdV_{i'}(\dul{a}^{m'})
{}.
\end{array}
\end{equation}
This lies in $\Fil_n \wom^1_X$ because $i'>n$ .
\end{proof}

As indicated in \cite{S2} \S 3.1 one can define for a morphism $f:X\rightarrow
S$ a relativized version of the De Rham-Witt complexes.
Here are the details for the \emph{relative
generalized De Rham-Witt complex of level $n$}.

Since the constructions are functorial $f$ induces a subsheaf
$f^*\won{n}_S\subset\won{n}_X$.
We let $\cJ^\bullet_{X/S,n}$ denote the graded ideal in
$\womn{n}^\bullet_X$ generated by
$d(f^*\won{n}_S)$; i.e.
\begin{equation}\label{eq:rel drw ideal}
\cJ^0_{X/S,n}=0\,,\qquad \cJ^i_{X/S,n}=
d(f^*\won{n}_S)\cdot\womn{n}^{i-1}_X\quad\textrm{for}\;i\geq 1,
\end{equation}
and define
\begin{equation}\label{eq:reldrw}
\womn{n}^\bullet_{X/S}:=
\womn{n}^\bullet_X/\cJ^\bullet_{X/S,n}.
\end{equation}
Note in particular
$$
\womn{n}^0_{X/S}=\won{n}_X.
$$
Clearly
$\womn{n}^\bullet_{X/S}$ is a differential graded algebra over
$\won{n}_S$.

For sections $\alpha$ of $f^*\won{n}_S$ and $\beta$ of
$\womn{n}^{i-1}_X$ one immediately computes
$$
V_m((d\alpha)\beta)=V_m((F_mdV_m\alpha)\beta)=(dV_m\alpha)(V_m\beta),
$$
while the more tricky calculation (\ref{eq:trick}) shows
$$
F_m(d(f^*\won{mn}_S)\cdot\womn{mn}^{i-1}_X)\subset
 d(f^*\won{n}_S)\cdot\womn{n}^{i-1}_X.
$$
Thus $V_m$ and $F_m$ induce maps
\begin{eqnarray*}
V_m:&&\womn{n}^i_{X/S}\rightarrow\womn{nm+m-1}^i_{X/S}
\\
F_m:&&\womn{mn}^i_{X/S}\rightarrow\womn{n}^i_{X/S}
\end{eqnarray*}
Exactly as in (\ref{eq:dr=drw1})
very mild conditions on $X$ will remove the $i!$-torsion from
$\OM^i_{X/S}$ and $\OM^i_{X}$, and imply that
$\womn{1}^\bullet_{X}\,=\,\OM^\bullet_{X}$,
$\;\cJ^\bullet_{X/S,1}=\cJ^\bullet_{X/S}$ (see (\ref{eq:JXS})) and
\begin{equation}\label{eq:reldr=drw1}
\womn{1}^\bullet_{X/S}\,=\,\OM^\bullet_{X/S}\,.
\end{equation}

\

\section{Characterization of $Z_N\OM^i_{X/S}$
by Frobenius and Verschiebung}\label{section:ZN frob}

Again we consider a morphism $f:X\rightarrow S$. For simplicity we assume
$S=\Spec A$ for some ring $A$. Then each $\womn{n}^i_{X/S}$ is a module over
the ring $\wan{n}$. For every $N$ we have the ring homomorphism
$$
F_N: \wan{nN}\rightarrow\wan{n}
$$
making $\wan{n}$ an algebra over $\wan{nN}$. We now \emph{define}
\begin{equation}\label{eq:F star}
F^*_N\womn{nN}^i_{X/S}\::=\:\wan{n}\otimes_{\wan{nN}}\womn{nN}^i_{X/S}.
\end{equation}
This is then a module over $\wan{n}$. The map
$$
F_N: \womn{nN}^i_{X/S}\rightarrow\womn{n}^i_{X/S}
$$
induces a \emph{linear map} of modules over $\wan{n}$
\begin{equation}\label{eq: F linear}
F_N: F^*_N\womn{nN}^i_{X/S}\rightarrow\womn{n}^i_{X/S}\,,\qquad
F_N(a\otimes\omega)=a\cdot F_N\omega.
\end{equation}

\

The Cartier isomorphism is crucial for a good theory of the De Rham(-Witt)
complex for schemes which are smooth over a scheme of characteristic $p>0$; cf.
\cite{K,IR}. The following result (\ref{eq:image FN}), the proof of which uses
in an essential way the Cartier isomorphism in characteristic $p$, is the
closest we can get to this in our setting.
For our purpose it gives one good characterization of $Z_N\OM^i_{X/S}$
in terms of the generalized De Rham-Witt complex.

\begin{proposition}\label{prop:Cartier iso}
Let $f:X\rightarrow S=\Spec A$ be a smooth projective morphism and assume that
the ring $A$ is \'etale over a polynomial ring $\ZZ\half[x_1,\ldots,x_r]$.
Then for all $i$
\begin{equation}\label{eq:wom1=om}
\womn{1}^i_{X/S}=\OM^i_{X/S}.
\end{equation}
Recall from (\ref{eq:ZN}) the definition of $Z_N\OM^i_{X/S}$.
Then for $N\geq 1$ and for all $i$
\begin{equation}\label{eq:image FN}
Z_N\OM^i_{X/S}=\mathrm{image}(F_N: F^*_N\womn{N}^i_{X/S}\rightarrow
\OM^i_{X/S}).
\end{equation}
\end{proposition}
\begin{proof}
The smoothness of $X$ over $S$ and flatness of $A$ over $\ZZ$ imply that there
is no $i!$-torsion in $\OM^i_{X/S}$ and hence that
$\womn{1}_{X/S}^i\,=\,\OM^i_{X/S}$ (cf. (\ref{eq:dr=drw1})).

The inclusion $\subset$ in (\ref{eq:image FN}) follows from the basic relation
$dF_N=NF_Nd$.
So let us concentrate on the other inclusion and first prove it in case $N$ is
a power of some prime number $p$. Write $S_p=S\times\Spec(\ZZ/p\ZZ)$ and
$X_p=X\times\Spec(\ZZ/p\ZZ)$. Then $X_p$ is smooth over $S_p$ and
$$
\OM^\bullet_{X_p/S_p}=\OM^\bullet_{X/S}\,/\,p\cdot\OM^\bullet_{X/S}.
$$
We are going make heavy use of \cite{K} thm. 7.2. As for notation
in loc. cit.:
$$
\OM^i_{X_p^{(p)}/S_p}:=(A/pA)\otimes_A \OM^i_{X/S}
$$
where $A/pA$ is considered as a module over $A$ via the map $a\mapsto a^p\bmod
p$.

Now take for some $\nu\geq 1$
\begin{equation}\label{eq:omega exact mod}
\omega\in\OM^i_{X/S}\quad \textit{such that}\quad d\omega\equiv 0\bmod p^\nu\,.
\end{equation}
This implies that $\omega\bmod p$ is a closed
form in $\OM^i_{X_p/S_p}$. The theorem on the Cartier isomorphism in
characteristic $p$ (see \cite{K} thm. 7.2) shows therefore
\begin{equation}\label{eq:omega cartier}
\omega\bmod p = C^{-1}\alpha+d\beta\qquad\textrm{with}\quad
\alpha\in \OM^i_{X_p^{(p)}/S_p}\,,\;
\beta\in\OM^{i-1}_{X_p/S_p}
\end{equation}
where $C^{-1}$ is the \emph{inverse Cartier operator}.\\
Choose $s_j\in A$, sections $a_{kj}$ of $\cO_X$ and section $\widetilde{\beta}$
of $\OM^{i-1}_{X/S}$ such that
$$
\beta=\widetilde{\beta}\bmod p\,,\qquad
\alpha=\widetilde{\alpha}\bmod p
\qquad\textrm{with}\quad
\widetilde{\alpha}=\sum_j s_j\otimes a_{0j}da_{1j}\cdot\ldots\cdot da_{ij}.
$$
Recall the following formulas from (\ref{eq:FVrel})
\begin{equation}
\label{eq:explicit frob}
\begin{array}{l}
F_pV_p=p\,,\qquad F_pdV_p=d\,,\\
F_p(\dul{a_0}d\dul{a_1}\cdot\ldots\cdot\dul{a_i})=
a_0^p(a_1\cdot\ldots\cdot a_i)^{p-1} da_1\cdot\ldots\cdot da_i\,.
\end{array}
\end{equation}
Comparing these with the formulas for
$C^{-1}$ in \cite{K} thm. 7.2 we see that (\ref{eq:omega cartier}) can be
rewritten as
\begin{eqnarray*}
&&
\omega=F_p\omega_1
\qquad\textrm{with}\\
&&
\omega_1=\sum_j s_j\otimes
\dul{a_{0j}}d\dul{a_{1j}}\cdot\ldots\cdot\dul{a_{ij}}
+1\otimes dV_p\widetilde{\beta}+1\otimes V_p\widetilde{\gamma}
\;\in\; F^*_p\womn{p}^i_{X/S},
\end{eqnarray*}
some $\widetilde{\gamma}\in\OM^i_{X/S}$.
If $\nu=1$ we are done. If $\nu>1$ we note that (\ref{eq:omega exact mod})
and $d\omega= p F_p d\omega_1$ imply
$$
\sum_j s_j\otimes (a_{0j}a_{1j}\cdot\ldots\cdot a_{ij})^{p-1}
da_{0j}da_{1j}\cdot\ldots\cdot da_{ij}\:+\: d\widetilde{\gamma}
=F_p d\omega_1\equiv 0\bmod p.
$$
Looking at the first term we see that this means
$$
C^{-1}d\alpha\:+\: d\widetilde{\gamma}=0.
$$
{}From the theorem on the Cartier isomorphism (\cite{K} thm. 7.2) we can now
conclude $d\alpha=0$ and thus as before
$$
\widetilde{\alpha}=F_p\widetilde{\alpha_1}\qquad\textrm{with}\quad
\widetilde{\alpha_1}\in F^*_p\womn{p}^i_{X/S}
$$
Lift $\widetilde{\alpha_1}$ to some $\widetilde{\widetilde{\alpha_1}}\in
F_p^*\womn{p^2}^i_{X/S}=\wan{p}\otimes_{\wan{p^2}}\womn{p^2}^i_{X/S}$
such that $\widetilde{\alpha_1}=(\tr\otimes\tr)
\widetilde{\widetilde{\alpha_1}}$. Then
$$
\omega=F_{p^2}\left(\widetilde{\widetilde{\alpha_1}}
+dV_{p^2}\widetilde{\beta}\right)+p\widetilde{\gamma}
{}.
$$
{}From this we see $d\widetilde{\gamma}=0\bmod p$ and thus
$$
\widetilde{\gamma}=F_p\widetilde{\gamma_1}\qquad\textrm{with}\quad
\widetilde{\gamma_1}\in F^*_p\womn{p}^{i-1}_{X/S}.
$$
 Take
$\widetilde{\widetilde{\gamma_1}}\in F_p^*\womn{p^2}^i_{X/S}$ such that
$F_p\widetilde{\widetilde{\gamma_1}}=p\widetilde{\gamma_1}$.
Then
$$
\omega=F_{p^2}\left(\widetilde{\widetilde{\alpha_1}}+
\widetilde{\widetilde{\gamma_1}}+dV_{p^2}\widetilde{\beta}
\right).
$$
If $\nu=2$ we are done. If $\nu>2$ we go on the same way untill finally
$$
\omega=F_{p^\nu}\omega_{p^\nu}\qquad\textrm{with}\quad
\omega_{p^\nu}\in F^*_{p^\nu}\womn{p^\nu}^i_{X/S}.
$$

Now take an arbitrary positive integer $N$ with prime decomposition
$N=\prod_{l}p_l^{\nu_l}$. Consider
$$
\omega\in\OM^i_{X/S}\quad \textit{such that}\quad d\omega\equiv 0\bmod N\,.
$$
By the previous prime power case we know that for every $l$
$$
\omega=F_{p_l^{\nu_l}}\eta_l\qquad\textrm{with}\quad
\eta_l\in F^*_{p_l^{\nu_l}}\womn{p_l^{\nu_l}}^i_{X/S}
$$
Take integers $c_l$ such that $\sum_l c_l p_l^{-\nu_l}N=1$.
Then
$$
\omega=\sum_l c_l\,F_{p_l^{\nu_l}}( p_l^{-\nu_l}N\,\eta_l)
=F_{N}\left(\sum_l c_l \,\widetilde{\eta_l}\right)
$$
with $\widetilde{\eta_l}\in F^*_N\womn{N}^i_{X/S}$ such that
$F_{p_l^{-\nu_l}N}\widetilde{\eta_l}=p_l^{-\nu_l}N\eta_l$.
\end{proof}

\

A second good characterization of $Z_N\OM^i_{X/S}$
in terms of the generalized De Rham-Witt complex is provided by (\ref{eq:ZVdV})
below. This result is an analogue of the exact sequence
\cite{IR} II(1.2.2):
$$
0\rightarrow W\OM^{i-1}\stackrel{(F^n,-F^nd)}{\longrightarrow}
W\OM^{i-1}\oplus W\OM^i\stackrel{(dV^n+V^n)}{\longrightarrow}
W\OM^i\rightarrow W_n\OM^i\rightarrow 0
$$
which is of fundamental importance for the theory of the De Rham-Witt complex
in characteristic $p$.

\begin{proposition}\label{prop:ker V+dV}
Let $f:X\rightarrow S=\Spec A$ be a smooth projective morphism and assume that
the ring $A$ is \'etale over a polynomial ring $\ZZ\half[x_1,\ldots,x_r]$.
Take $N>1$ and recall from (\ref{eq:ZN}) the definition of $Z_N\OM^i_{X/S}$.
Then one has the exact sequences
\begin{equation}\label{eq:ZVdV0}
 0\rightarrow \OO_X\stackrel{V_N}{\longrightarrow}
\won{N}_X\longrightarrow\won{N-1}\rightarrow 0
\end{equation}
and for all $i\geq 1$
\begin{equation}\label{eq:ZVdV}
0\rightarrow Z_N\OM^{i-1}_{X/S}\stackrel{(1,-\frac{1}{N}d)}{\longrightarrow}
\OM^{i-1}_{X/S}\oplus \OM^i_{X/S}\stackrel{dV_N+V_N}{\longrightarrow}
\womn{N}^i_{X/S}\longrightarrow\womn{N-1}^i_{X/S}\rightarrow 0.
\end{equation}
with a slight abuse of notation, in that the map $V_N$ is actually
$\tr\circ V_N$ with
$\tr:\womn{2N-1}^\bullet_{X/S}\rightarrow\womn{N}^\bullet_{X/S}$ the standard
truncation.
\end{proposition}
\begin{proof}
It is immediately obvious from the definitions and constructions in Section
\ref{section:FrobWitt} that the sequence (\ref{eq:ZVdV0}) is exact.
The only point where it is not immediately obvious from the definitions and
constructions in Section \ref{section:DRW} that the sequence (\ref{eq:ZVdV}) is
exact, is at $\OM^{i-1}_{X/S}\oplus \OM^i_{X/S}$.

The composite map
$
\OM^{i-1}_{X/S}\oplus \OM^i_{X/S}\stackrel{dV_N+V_N}{\longrightarrow}
\womn{N}^i_{X/S}\stackrel{F_N}{\longrightarrow}\OM^i_{X/S}
$
is in fact $d+N$. This shows
$$
\ker\left(\OM^{i-1}_{X/S}\oplus \OM^i_{X/S}\stackrel{dV_N+V_N}{\longrightarrow}
\womn{N}^i_{X/S}\right)\subset
\mathrm{image}\left(
Z_N\OM^{i-1}_{X/S}\stackrel{(1,-\frac{1}{N}d)}{\longrightarrow}
\OM^{i-1}_{X/S}\oplus \OM^i_{X/S}\right)
$$
For the $\supset$-inclusion we note that by (\ref{eq:image FN}) with $i-1$ in
place of $i$ it suffices to show that the composite map
$$
F^*_N\womn{N}^{i-1}_{X/S}\stackrel{(F_N,-F_Nd)}{\longrightarrow}
\OM^{i-1}_{X/S}\oplus \OM^i_{X/S}\stackrel{dV_N+V_N}{\longrightarrow}
\womn{N}^i_{X/S}
$$
is zero. Since by definition
$F^*_N\womn{N}^{i-1}_{X/S}=A\otimes_{\wan{N}}\womn{N}^{i-1}_{X/S}$ this means
that we must show for $a\in A$ and $\omega\in \womn{N}^{i-1}_{X/S}$
$$
dV_N(aF_N\omega)=V_N(aF_Nd\omega).
$$
With the relations in (\ref{eq:FVrel}) the left hand side can  be written as
$d((V_Na)\omega)$, while the right hand  equals $(V_Na)d\omega$.
These two are equal since $V_Na\in\wan{N}$. This completes the proof.
\end{proof}

With the techniques of the preceding proofs we can now also prove the claim
made at the end of Remark \ref{remark:ordinary open}.

\begin{proposition}\label{prop:reduce ordinary}
Let $f:X\rightarrow S=\Spec A$ be a smooth projective morphism such that the
ring $A$ is \'etale over a polynomial ring $\ZZ\half[x_1,\ldots,x_r]$ and
$H_{DR}^m(X/S)$ and
$H^{m-i}(X,\OM_{X/S}^i)$ are free $A$-modules for all $m,i$.
Fix $N\geq 2$ and
assume that
\begin{equation}\label{eq:p}
\begin{array}{rl}
H^j(X,Z_p\OM_{X/S}^i)\rightarrow H^j(X,\OM_{X/S}^i)
&\textit{is an isomorphism for all }\; i,j\;
\textit{ and}\\
&\textit{all prime numbers }\;p\;\textit{ dividing }\; N.
\end{array}
\end{equation}
Then
\begin{equation}\label{eq:N}H^j(X,Z_N\OM_{X/S}^i)\rightarrow H^j(X,\OM_{X/S}^i)
\quad\textrm{is an isomorphism for all}\quad i,j.
\end{equation}

\end{proposition}
\begin{proof}\textbf{Step 1.}
Consider the prime decomposition $N=\prod_l p_l^{\nu_l}$. Fix $i,j$.
Assume that for every $p_l$ with $\nu_l>0$
the map $H^j(X,Z_{p_l^{\nu_l}}\OM_{X/S}^i)\rightarrow
H^j(X,\OM_{X/S}^i)$ is an isomorphism and call this map $\rho_l$.
Note that multiplication by $Np_l^{-\nu_l}$ induces a map
$$
H^j(X,Z_{p_l^{\nu_l}}\OM_{X/S}^i)
\stackrel{\cdot Np_l^{-\nu_l}}{\longrightarrow}
H^j(X,Z_N\OM_{X/S}^i).
$$
Choose integers $c_l$ such that $\sum_l c_lNp_l^{-\nu_l}\,=\,1$.
Then
$$
\sum_l c_lNp_l^{-\nu_l}\rho_l^{-1}:H^j(X,\OM_{X/S}^i)\rightarrow
H^j(X,Z_N\OM_{X/S}^i)
$$
is an inverse for the map in (\ref{eq:N}).

\textbf{Step 2.} The problem is thus reduced to showing that if $p$ is a prime
number for which (\ref{eq:p}) holds, then for every $\nu\geq 1$
\begin{equation}\label{eq:pnu}
H^j(X,Z_{p^\nu}\OM_{X/S}^i)\rightarrow H^j(X,\OM_{X/S}^i)
\qquad\textrm{is an isomorphism for all}\quad i,j.
\end{equation}
The basic relationship between the Frobenius operator $F_p$ and the inverse
Cartier operator $C^{-1}$ is expressed in the following commutative diagram
\begin{equation}\label{eq:frob-cartier}
\begin{array}{ccc}
F^*_p\womn{p}^i_{X/S}&\stackrel{F_p}{\longrightarrow}&Z_p\OM^i_{X/S} \\
\downarrow&&\downarrow \\
(A/pA)\otimes_A \OM^i_{X/S}&\stackrel{C^{-1}}{\longrightarrow}&
(Z_p\OM^i_{X/S})/(d\OM^{i-1}_{X/S}+p\OM^i_{X/S})
\end{array}
\end{equation}
where the left hand vertical map is the one induced by the standard truncation
$\womn{p}^i_{X/S}\rightarrow\OM^i_{X/S}$ and the right hand vertical map is the
obvious one. In the notation $(A/pA)\otimes_A$ the ring $A/pA$ is an $A$-module
via the ring homomorphism
$A\rightarrow A/pA\,,\;a\mapsto a^p\bmod p$.
Composing with $F_{p^\nu}:\:
F_{p^{\nu+1}}^*\womn{p^{\nu+1}}\rightarrow F^*_p\womn{p}^i_{X/S}$ in the upper
left hand corner we obtain the commutative diagram
$$
\begin{array}{ccc}
F_{p^{\nu+1}}^*\womn{p^{\nu+1}}&\stackrel{=}{\longrightarrow}&
F_{p^{\nu+1}}^*\womn{p^{\nu+1}}
\\[1.5ex]
\hspace{-2em}F_{p^\nu}\downarrow&&\hspace{2.5em}\downarrow F_{p^{\nu+1}}\\
F^*_p\womn{p}^i_{X/S}&\stackrel{F_p}{\longrightarrow}&Z_p\OM^i_{X/S} \\
\downarrow&&\downarrow \\
(A/pA)\otimes_A \OM^i_{X/S}&\stackrel{C^{-1}}{\longrightarrow}&
(Z_p\OM^i_{X/S})/(d\OM^{i-1}_{X/S}+p\OM^i_{X/S})
\end{array}
$$
According to \cite{K} thm. 7.2 the inverse Cartier operator $C^{-1}$
is an isomorphism.
Taking cokernels of the vertical composite maps and using
(\ref{eq:image FN}) one obtains an isomorphism
$$
(A/pA)\otimes_A
\begin{array}{c}
\OM^i_{X/S}
\\
\overline{Z_{p^{\nu}}\OM^i_{X/S}+p\OM^i_{X/S}}
\end{array}
\quad\stackrel{\simeq}{\longrightarrow}\quad
\begin{array}{c}
Z_p\OM^i_{X/S}
\\
\overline{Z_{p^{\nu+1}}\OM^i_{X/S}+p\OM^i_{X/S}}
\end{array}\,.
$$
On cohomology this gives an isomorphism
$$
H^j\left(X,
(A/pA)\otimes_A
\begin{array}{c}
\OM^i_{X/S}
\\
\overline{Z_{p^{\nu}}\OM^i_{X/S}+p\OM^i_{X/S}}
\end{array}\right)
\stackrel{\simeq}{\longrightarrow}
H^j\left(X,\begin{array}{c}
Z_p\OM^i_{X/S}
\\
\overline{Z_{p^{\nu+1}}\OM^i_{X/S}+p\OM^i_{X/S}}
\end{array}\right)\,.
$$
Since the Frobenius endomorphism $A/pA\rightarrow A/pA\,,\;a\mapsto a^p$,
is a flat map the group on the left hand side is isomorphic to
$$
(A/pA)\otimes_A
H^j\left(X,
\begin{array}{c}
\OM^i_{X/S}
\\
\overline{Z_{p^{\nu}}\OM^i_{X/S}+p\OM^i_{X/S}}
\end{array}\right).
$$
Hypothesis (\ref{eq:p}) implies that the group on the right hand side is
isomorphic to
$$
H^j\left(X,\begin{array}{c}
\OM^i_{X/S}
\\
\overline{Z_{p^{\nu+1}}\OM^i_{X/S}+p\OM^i_{X/S}}
\end{array}\right).
$$
These arguments show in particular:
$$
\textit{if}\;
H^j\left(X,
\begin{array}{c}
\OM^i_{X/S}
\\
\overline{Z_{p^{\nu}}\OM^i_{X/S}+p\OM^i_{X/S}}
\end{array}\right)=0\;\,
\textit{then}\;
H^j\left(X,\begin{array}{c}
\OM^i_{X/S}
\\
\overline{Z_{p^{\nu+1}}\OM^i_{X/S}+p\OM^i_{X/S}}
\end{array}\right)=0.
$$
This implication can be combined with implications resulting from the exact
cohomology sequences of the short exact sequences of sheaves
$$
0\rightarrow
\begin{array}{c}
\OM^i_{X/S}
\\
\overline{Z_{p^r}\OM^i_{X/S}}
\end{array}
\stackrel{\cdot p}{\longrightarrow}
\begin{array}{c}
\OM^i_{X/S}
\\
\overline{Z_{p^{r+1}}\OM^i_{X/S}}
\end{array}
\longrightarrow
\begin{array}{c}
\OM^i_{X/S}
\\
\overline{Z_{p^{r+1}}\OM^i_{X/S}+p\OM^i_{X/S}}
\end{array}
\rightarrow
0
$$
for every $r\geq 0$. This leads to the conclusion:
\begin{eqnarray*}
\textit{If}&\quad&
H^j\left(X,\begin{array}{c}
\OM^i_{X/S}
\\
\overline{Z_{p^{\nu-1}}\OM^i_{X/S}}
\end{array}
\right)=
H^j\left(X,\begin{array}{c}
\OM^i_{X/S}
\\
\overline{Z_{p^\nu}\OM^i_{X/S}}
\end{array}
\right)=0\qquad\textit{ for all }\;j,\\
\textit{then}&\quad&
H^j\left(X,\begin{array}{c}
\OM^i_{X/S}
\\
\overline{Z_{p^{\nu+1}}\OM^i_{X/S}}
\end{array}
\right)=0\qquad\textit{ for all }\;j.
\end{eqnarray*}
Thus we can derive (\ref{eq:pnu}) from (\ref{eq:p}) by induction.
\end{proof}

\

\section{The conjugate filtration and ordinariness: Act 2}
\label{section:conjugate2}
\subsection{The generalized Hodge-Witt complex.}
\label{subsection:Hodge-Witt}
Let $X$ be a scheme on which $2$ is invertible.
For $n\in\NN$ we define the \emph{generalized Hodge-Witt complex of level $n$
on $X$} to be the graded algebra
\begin{equation}\label{eq:hodgewittcomplex}
\womn{n}^\oplus_X:=
\bigoplus_{i\geq 0} \womn{n}_X^i[-i]
\end{equation}
with zero differential. As graded algebras $\womn{n}^\oplus_X$
and $\womn{n}^\bullet_X$ are the same, but they carry different differentials.
The fundamental relation
$$
dF_N=NF_Nd
$$
implies that the maps $F_N:\womn{nN}_X^i\rightarrow \womn{n}_X^i$
together yield, for every $n$, a
\emph{homorphism of differential graded algebras}
\begin{equation}\label{eq:basic morphism 1}
\wtF_N:\;\womn{nN}^\oplus_X\;\longrightarrow\;
\womn{n}^\bullet_X\q{N}.
\end{equation}

Not only is $\wtF_N$  a homorphism of differential graded algebras, but
it preserves several other important algebraic structures: see
(\ref{eq:tF trunc})--(\ref{eq:basic morphism 2}).
\\
For $k>n$ one has the commutative square
\begin{equation}\label{eq:tF trunc}
\begin{array}{ccc}
\womn{kN}^\oplus_X&\stackrel{\wtF_N}{\longrightarrow}&
\womn{k}^\bullet_X\q{N}\\[1ex]
\tr\downarrow&&\downarrow\tr\\[1ex]
\womn{nN}^\oplus_X&\stackrel{\wtF_N}{\longrightarrow}&
\womn{n}^\bullet_X\q{N}
\end{array}
\end{equation}
in which the vertical maps are the canonical truncation maps.
\\
For $K|N$ one has the commutative square
\begin{equation}\label{eq:tF F}
\begin{array}{ccc}
\womn{nN}^\oplus_X&\stackrel{\wtF_N}{\longrightarrow}&
\womn{n}^\bullet_X\q{N}\\[1ex]
\hspace{-2em}F_{N/K}\downarrow&&\downarrow\tr\\[1ex]
\womn{nK}^\oplus_X&\stackrel{\wtF_K}{\longrightarrow}&
\womn{n}^\bullet_X\q{K}
\end{array}
\end{equation}
\\
For $K|N$ one has the commutative diagram
\begin{equation}\label{eq:tF V}
\begin{array}{ccc}
\womn{nK}^\oplus_X&\stackrel{\wtF_K}{\longrightarrow}&
\womn{n}^\bullet_X\q{K}\\[1ex]
\hspace{-4em}\tr\circ V_{N/K}\downarrow&&\hspace{1em}\downarrow \cdot
N/K\\[1ex]
\womn{nN}^\oplus_X&\stackrel{\wtF_N}{\longrightarrow}&
\womn{n}^\bullet_X\q{N}
\end{array}
\end{equation}
\\
The fundamental relation
$$
dF_m=mF_md
$$
also shows that there is a
\emph{homomorphism of differential graded algebras}
\begin{equation}\label{eq:FF}
\FF_m:\womn{nm}^\bullet_X\longrightarrow \womn{n}^\bullet_X\,,
\qquad \FF_m\alpha=m^iF_m\alpha\quad\textrm {for}\quad\alpha\in
\womn{nm}^i_X
\end{equation}
and similarly for $\womn{.}^\oplus_X$ in place of $\womn{.}^\bullet_X$.
The following square is commutative
\begin{equation}\label{eq:tF frob}
\begin{array}{ccc}
\womn{mnN}^\oplus_X&\stackrel{\wtF_N}{\longrightarrow}&
\womn{mn}^\bullet_X\q{N}\\[1ex]
\hspace{-1em}\FF_m\downarrow&&\downarrow\FF_m\\[1ex]
\womn{nN}^\oplus_X&\stackrel{\wtF_N}{\longrightarrow}&
\womn{n}^\bullet_X\q{N}
\end{array}
\end{equation}
\\
Recall that we defined the ideal $\cJ^\bullet_{X/S,n}$ in $\womn{n}^\bullet_X$
by the formulas in (\ref{eq:rel drw ideal}). With the same formulas we define
the ideal $\cJ^\oplus_{X/S,n}$ in $\womn{n}^\oplus_X$. So as ideals
$\cJ^\oplus_{X/S,n}$ and $\cJ^\bullet_{X/S,n}$ are the same, but they carry
different differentials.
The map $\wtF_N$ restricts to a homomorphism of differential graded ideals
\begin{equation}\label{eq:J HW DRW}
\wtF_N:\;\cJ^\oplus_{X/S,nN}\longrightarrow\cJ^\bullet_{X/S,n}\q{N}.
\end{equation}
Modding out by these ideals we get the \emph{relative generalized Hodge-Witt
(resp. De Rham-Witt) complex of level $n$}:
\begin{equation}\label{eq:relative}
\womn{n}^\oplus_{X/S}:=\womn{n}^\oplus_X/\cJ^\oplus_{X/S,n}\:,\qquad
\womn{n}^\bullet_{X/S}:=\womn{n}^\bullet_X/\cJ^\bullet_{X/S,n}.
\end{equation}
The map $\wtF$ induces a homomorphism of differential graded algebras
\begin{equation}\label{eq:basic morphism 2}
\wtF_N:\;\womn{nN}^\oplus_{X/S}\;\longrightarrow\;
\womn{n}^\bullet_{X/S}\q{N}.
\end{equation}
\textbf{Fact:}
There are obvious analogues of (\ref{eq:tF trunc})--(\ref{eq:tF frob}) with
$X/S$ instead of $X$.

\begin{remark}
If $X$ is a smooth scheme over a perfect field of characteristic $p>2$, if
$n=N=p^j$ and if we restrict to the $p$-typical parts, then $\wtF_N$ induces on
the cohomology sheaves the \emph{higher Cartier isomorphism} $C^{-j}$ of
\cite{IR} p. 77.
In the same situation, with $j=1$, the map $\wtF_p$ plays a decisive role in
Deligne and Illusie's algebraic proof of the degeneration of the Hodge-to-De
Rham spectral sequence; see \cite{DI} thm. 2.1.
\end{remark}

\subsection{From Hodge-Witt to De Rham-Witt.}
\label{subsection:split hodge}

{}From now on we assume that
$f:X\rightarrow S=\Spec A$ is a smooth projective morphism in which the ring
$A$ is \'etale over a polynomial ring $\ZZ\half[x_1,\ldots,x_r]$ and that
$H_{DR}^m(X/S)$ and
$H^{m-i}(X,\OM_{X/S}^i)$ are free $A$-modules for all $m,i$.

Note that the smoothness of $X$ over $S$ and flatness of $A$ over $\ZZ$ imply
that there is no $i!$-torsion in $\OM^i_{X/S}$ and hence that
$$
\womn{1}_{X/S}^\bullet\,=\,\OM^\bullet_{X/S}.
$$
Fix $N\in \NN$. The homomorphism $\wtF_N$ from (\ref{eq:basic morphism 2})
induces on the $m$-th hypercohomology the map
\begin{equation}\label{eq:basic morphism 3}
\wtF_N: \bigoplus_j H^{m-j}(X,\womn{N}^j_{X/S})\longrightarrow
\HH^m(X,\OM^\bullet_{X/S}\q{N}).
\end{equation}
Since the De Rham cohomology groups are free $A$-modules and  $A$ is flat over
$\ZZ$ we have
$$
\HH^m(X,\OM^\bullet_{X/S}\q{N})=H^m_{DR}(X/S)\q{N}.
$$

\begin{proposition}\label{prop:PhiN compatibilities}
$\Phi_N$ relates as follows to the structures on $H^m_{DR}(X/S)\q{N}$ discussed
in Section \ref{section:conjugate}.
\\
$\bullet$
\textup{Hodge filtration:   }
For all $N,m,i$
\begin{equation}\label{eq:sub Hodge filt}
\wtF_N\left( \bigoplus_{j\geq i} H^{m-j}(X,\womn{N}^j_{X/S})\right)
\subset \FH^i H^m_{DR}(X/S)\q{N}.
\end{equation}
\\
$\bullet$
\textup{Conjugate filtration:   }
For all $N,m,i$
\begin{equation}\label{eq:sub con filt}
\wtF_N\left( \bigoplus_{j\leq i} H^{m-j}(X,\womn{N}^j_{X/S})\right)
\subset \FC_i\HH^m(X,\OM^\bullet_{X/S}\q{N}).
\end{equation}
\\
$\bullet$
\textup{Gauss-Manin connection:   }
For all $N,j$
\begin{equation}\label{eq: G-M nilpotent}
\nabla\circ\wtF_N\left(H^{m-j}(X,\womn{N}^j_{X/S})\right)\;\subset
\; \OM_{S/\ZZ}^1\otimes\wtF_N\left(H^{m-j+1}(X,\womn{N}^{j-1}_{X/S})
\right)
\end{equation}
\end{proposition}
\begin{proof}
(\ref{eq:sub Hodge filt}) and (\ref{eq:sub con filt}) are obvious.
Let us prove (\ref{eq: G-M nilpotent}).
Recall from (\ref{eq:J HW DRW}) that $\wtF_N$ maps $\cJ^\oplus_{X/S,N}$ into
$\cJ^\bullet_{X/S,1}\q{N}$. So there is a commutative diagram with
exact rows and vertical maps $\wtF_N$:
$$
\begin{array}{ccccccccc}
0&\!\!\rightarrow\!\!&\cJ^\oplus_{X/S,N}&\!\!\rightarrow\!\!&
\womn{N}^\oplus_X
&\!\!\rightarrow\!\!&\womn{N}^\oplus_{X/S}&\!\!\rightarrow\!\!&0\\
&&\downarrow&&\downarrow&&\downarrow&&\\
0&\!\!\rightarrow\!\!& \cJ^\bullet_{X/S,1}/(\cJ^\bullet_{X/S,1})^2\q{N}
&\!\!\rightarrow\!\!
&\OM_X^\bullet/(\cJ^\bullet_{X/S,1})^2\q{N}&\!\!\rightarrow\!\!&
\OM_{X/S}^\bullet\q{N}&\rightarrow &0
\end{array}
$$
In the ladder of hypercohomology groups there is the commutative square
$$
\begin{array}{ccc}
\HH^m(X,\womn{N}^\oplus_{X/S})&\rightarrow&
\HH^{m+1}(X,\cJ^\oplus_{X/S,N})\\
\downarrow&&\downarrow\\
\HH^m(X,\OM_{X/S}^\bullet\q{N})&\rightarrow&
\HH^{m+1}(X,\cJ^\bullet_{X/S,1}/(\cJ^\bullet_{X/S,1})^2\q{N})
\end{array}
$$
in which the bottom map is in fact the Gauss-Manin connection.
To finish the proof we note that $\cJ^0_{X/S,N}=0$ while for every $j\geq 1$
\begin{eqnarray*}
\wtF_N(\cJ^j_{X/S,N})&=&\wtF_N(d(f^*\won{N}_S)\cdot\womn{N}^{j-1}_{X})\\
&=&\wtF_N(d(f^*\won{N}_S))\cdot\wtF_N(\womn{N}^{j-1}_{X})
\subset\OM^1_S\cdot\wtF_N(\womn{N}^{j-1}_{X}).
\end{eqnarray*}
\end{proof}

\

\

Each $H^j(X,\womn{N}^i_{X/S})$ is a module over
the ring $\wan{N}$. Viewing $A$ as an algebra over $\wan{N}$ via the ring
homomorphism
$$
F_N: \wan{N}\rightarrow A
$$
we define
\begin{equation}\label{eq:F star H}
F^*_N H^j(X,\womn{N}^i_{X/S})\::=\:A\otimes_{\wan{N}}
H^j(X,\womn{N}^i_{X/S}).
\end{equation}
This is a module over $A$. The map
$\wtF_N$ in (\ref{eq:basic morphism 3})
induces a \emph{$A$-linear map}
\begin{eqnarray}\label{eq: F linear2}
&&\wtF_N: \bigoplus_j F^*_N H^{m-j}(X,\womn{N}^j_{X/S})\longrightarrow
H^m_{DR}(X/S)\q{N}\\
\nonumber
&&\hspace{5em}
\wtF_N(a\otimes\omega)=a\cdot \wtF_N\omega.
\end{eqnarray}

\

\begin{theorem}\label{thm:U=conj}
Let $f:X\rightarrow S=\Spec A$ be a smooth projective morphism with ring $A$
\'etale over a polynomial ring $\ZZ\half[x_1,\ldots,x_r]$.
Assume that $H_{DR}^m(X/S)$ and
$H^{m-i}(X,\OM_{X/S}^i)$ are free $A$-modules for all $m,i$. Fix $N\in\NN$ and
assume that
$H^j(X,Z_N\OM_{X/S}^i)\rightarrow H^j(X,\OM_{X/S}^i)$
is an isomorphism for all $ i,j$.
Fix $m$. Then the following three statements are equivalent
\\
\textup{(a)}\hspace{2em}
For $i=0,\ldots,m$:
\begin{equation}\label{eq:is Hodge filt}
\wtF_N\left( \bigoplus_{j\geq i} F^*_NH^{m-j}(X,\womn{N}^j_{X/S})\right)
= \FH^i H^m_{DR}(X/S)\q{N}.
\end{equation}
\textup{(b)}\hspace{2em}For $i=0,\ldots,m$:
\begin{equation}\label{eq:is con filt}
\wtF_N\left( \bigoplus_{j\leq i} F^*_NH^{m-j}(X,\womn{N}^j_{X/S})\right)
= \FC_i\HH^m(X,\OM^\bullet_{X/S}\q{N}).
\end{equation}
\textup{(c)}\hspace{2em}For $i=0,\ldots,m$ the map
$$
F_N: F^*_NH^{m-i}(X,\womn{N}^i_{X/S})\rightarrow
H^{m-i}(X,\OM^i_{X/S})\qquad\textrm{is surjective.}
$$
\end{theorem}
\begin{proof}
{}From the discussion preceding Theorem \ref{thm:conj split} we know that
the hypotheses imply for $i=0,\ldots,m$
\begin{eqnarray*}
\FH^i\HH^m(X,\OM^\bullet_{X/S}\q{N})/
\FH^{i+1}\HH^m(X,\OM^\bullet_{X/S}\q{N})
&=&H^{m-i}(X,\OM^i_{X/S})\q{N},\\
\FC_i\HH^m(X,\OM^\bullet_{X/S}\q{N})/
\FC_{i-1}\HH^m(X,\OM^\bullet_{X/S}\q{N})
&=&H^{m-i}(X,\OM^i_{X/S})\q{N}.
\end{eqnarray*}
On the other hand it is obvious that $\wtF_N$ induces maps
\begin{eqnarray*}
F^*_NH^{m-i}(X,\womn{N}^i_{X/S})&\hspace{-1em}\rightarrow&\hspace{-1em}
\FH^i\HH^m(X,\OM^\bullet_{X/S}\q{N})/
\FH^{i+1}\HH^m(X,\OM^\bullet_{X/S}\q{N}),
\\
F^*_NH^{m-i}(X,\womn{N}^i_{X/S})&\hspace{-1em}\rightarrow&\hspace{-1em}
\FC_i\HH^m(X,\OM^\bullet_{X/S}\q{N})/
\FC_{i-1}\HH^m(X,\OM^\bullet_{X/S}\q{N}).
\end{eqnarray*}
These are in fact equal to the map
\begin{equation}\label{eq:FN surj}
F^*_NH^{m-i}(X,\womn{N}^i_{X/S})\rightarrow H^{m-i}(X,\OM^i_{X/S})\q{N}
\end{equation}
which is $F_N$ followed by reduction $\bmod N$.
Thus we see that statements (a) and (b) are both equivalent to the statement
that the map in (\ref{eq:FN surj}) is surjective for $i=0,\ldots,m$.
The latter statement is clearly also implied by (c).
Conversely, since
$$
H^{m-i}(X,\OM^i_{X/S})\stackrel{V_N}{\longrightarrow}
H^{m-i}(X,\womn{N}^i_{X/S})\stackrel{F_N}{\longrightarrow}
H^{m-i}(X,\OM^i_{X/S})
$$
is just multiplication by $N$, we see that $NH^{m-i}(X,\OM^i_{X/S})$ is
contained in the image of $F_N$ and thus we can conclude that if the map in
(\ref{eq:FN surj}) is surjective, then
$$
F_N:\:F^*_NH^{m-i}(X,\womn{N}^i_{X/S})\rightarrow H^{m-i}(X,\OM^i_{X/S})
$$
is surjective too, i.e. (c) holds.
\end{proof}

\

\subsection{Compatibilities as $N$ varies.}\label{proj-ind systems}

{}From (\ref{eq:tF trunc}) we get for $K|N$ a commutative square
\begin{equation}\label{eq:prosystems}
\begin{array}{ccc}
\bigoplus_j  H^{m-j}(X,\womn{N}^j_{X/S})
&\stackrel{\wtF_N}{\longrightarrow}&
H^m_{DR}(X/S)\q{N}\\[1ex]
F_{N/K}\downarrow&&\downarrow\bmod K \\[1ex]
\bigoplus_j  H^{m-j}(X,\womn{K}^j_{X/S})
&\stackrel{\wtF_K}{\longrightarrow}&
H^m_{DR}(X/S)\q{K}
\end{array}
\end{equation}
On one side we have the projective system of groups
$\{H^m_{DR}(X/S)\q{N}\}_{N\in\NN}\}$ with for $K|N$ the reduce-$\bmod K$-map.
On the other side
we must consider the projective system of groups
$\{H^{m-j}(X,\womn{N}^j_{X/S})\}_{N\in\NN}$ with for $K|N$ the map $F_{N/K}$.

{}From (\ref{eq:tF V}) we get for $K|N$ a commutative square
\begin{equation}\label{eq:indsystems}
\begin{array}{ccc}
\bigoplus_j  H^{m-j}(X,\womn{N}^j_{X/S})
&\stackrel{\wtF_N}{\longrightarrow}&
H^m_{DR}(X/S)\q{N}\\[1ex]
V_{N/K}\uparrow&&\uparrow \cdot N/K \\[1ex]
\bigoplus_j  H^{m-j}(X,\womn{K}^j_{X/S})
&\stackrel{\wtF_K}{\longrightarrow}&
H^m_{DR}(X/S)\q{K}
\end{array}
\end{equation}
On one side we have the inductive system of groups
$\{H^m_{DR}(X/S)\q{N}\}_{N\in\NN}$ with for $K|N$
the map induced by multiplication by $N/K$.
On the other side we have the inductive system of groups
$\{H^{m-j}(X,\womn{N}^j_{X/S})\}_{N\in\NN}$ with for $K|N$ the map $V_{N/K}$
(or rather $\tr^{N+N/K-1}_N\circ V_{N/K}$).

In connection with the systems on the right it may be worthwhile to note that
$$
\ZZ/N\ZZ\;\simeq \frac{1}{N}\ZZ/\ZZ\,,\qquad z\bmod N\ZZ\:\leftrightarrow\:
\frac{z}{N}\bmod \ZZ
$$
and that via this isomorphism the map
$\bmod K:\ZZ/N\ZZ\rightarrow\ZZ/K\ZZ$ corresponds with the multiplication by
$N/K:$ $\frac{1}{N}\ZZ/\ZZ\rightarrow\frac{1}{K}\ZZ/\ZZ$
and the map
$\cdot N/K:\ZZ/K\ZZ\rightarrow\ZZ/N\ZZ$ corresponds with the inclusion
$\frac{1}{K}\ZZ/\ZZ\subset\frac{1}{N}\ZZ/\ZZ$.
\

\section{The conjugate filtration and ordinariness: Act 3}
\label{section:conjugate3}

\subsection{}

\

The general hypotheses for this section are:
$f:X\rightarrow S=\Spec A$ is a smooth projective morphism such that the ring
$A$ is \'etale over a polynomial ring $\ZZ\half[x_1,\ldots,x_r]$ and
$H_{DR}^m(X/S)$ and
$H^{m-i}(X,\OM_{X/S}^i)$ are free $A$-modules for all $m,i$.

We are going to analyse the projective systems of groups
$\{H^j(X,\womn{n}^i_{X/S})\}_{n\in\NN}$ with for $k<n$ the map induced by the
truncation $\tr^n_k:\womn{n}^i_{X/S}\rightarrow\womn{k}^i_{X/S}$. Later we look
at the issue of surjectivity of the maps
$$
F_N: F^*_NH^j(X,\womn{N}^i_{X/S})\rightarrow H^j(X,\OM^i_{X/S}).
$$

\

\begin{proposition}\label{prop:wittsurjective1}
Under the general hypotheses of this section the sequence
\begin{equation}\label{eq:V0}
0\rightarrow H^j(X,\OO_X)\stackrel{V_N}{\longrightarrow}
  H^j(X,\won{N}_X)\stackrel{\tr}{\longrightarrow} H^j(X,\won{N-1}_X)\rightarrow
0
\end{equation}
is exact and the map, induced by truncation,
\begin{equation}\label{eq:surj N0}
H^j(X,\womn{N}^i_{X/S})\rightarrow H^j(X,\OM^i_{X/S})
\end{equation}
is surjective for all $j$ and all $N\geq 2$.
\\
\textup{Note: we make here a slight abuse of notation by writing $V_N$ instead
of $\tr^{2N-1}_N\circ V_N$.}
\end{proposition}
\begin{proof}
Take $N\geq 2$. Consider the exact sequence
$$
 0\rightarrow \OO_X\stackrel{V_N}{\longrightarrow}
\won{N}_X\stackrel{\tr}{\longrightarrow}\won{N-1}\rightarrow 0
$$
and the following piece of the associated sequence of cohomology groups
$$
\rightarrow H^j(X,\OO_X)\stackrel{V_N}{\longrightarrow}
  H^j(X,\won{N}_X)\longrightarrow H^j(X,\won{N-1}_X)\rightarrow
  H^{j+1}(X,\OO_X)\stackrel{V_N}{\longrightarrow}
$$
One has for every $k$ the map $F_N:H^k(X,\won{N}_X)\rightarrow H^k(X, \OO_X)$
and $F_NV_N = N$. Since $H^k(X,\OO_X)$ is a free $A$-module and $A$ is flat
over $\ZZ$,  multiplication by $N$ is injective. Hence in the exact sequence
of cohomology groups all maps $V_N$ are injective and the long exact cohomology
sequence splits up into the short exact sequences
(\ref{eq:V0}).
\end{proof}

\begin{proposition}\label{prop:wittsurjective2}
Fix $N\geq 2$. Assume in addition to the general hypotheses of this section
that $H^j(X,Z_N\OM_{X/S}^i)\rightarrow H^j(X,\OM_{X/S}^i)$
is an isomorphism for all $i,j$.
\\
Then the sequence
\begin{equation}\label{eq:Vi}
0\rightarrow H^j(X,\OM^i_{X/S})\stackrel{V_N}{\longrightarrow}
H^j(X,\womn{N}^i_{X/S})\rightarrow H^j(X,\womn{N-1}^i_{X/S})\rightarrow 0
\end{equation}
is exact for all $i,j$.
\end{proposition}
\begin{proof}
Fix $i$. Set $\KK:=\ker(\womn{N}^i_{X/S}\rightarrow\womn{N-1}^i_{X/S})$. The we
have according to Proposition \ref{prop:ker V+dV} exact sequences
\begin{eqnarray}
\label{eq:right K}
&&
0\rightarrow Z_N\OM^{i-1}_{X/S}\stackrel{(1,-\frac{1}{N}d)}{\longrightarrow}
\OM^{i-1}_{X/S}\oplus
\OM^i_{X/S}\stackrel{dV_N+V_N}{\longrightarrow}\KK\rightarrow 0\\
\label{eq:left K}
&&
0\rightarrow\KK\rightarrow
\womn{N}^i_{X/S}\longrightarrow\womn{N-1}^i_{X/S}\rightarrow 0
\end{eqnarray}
The additional hypothesis and the exact cohomology sequence of
(\ref{eq:right K}) together show that the maps
\begin{equation}\label{eq:ker K}
H^j(X,\OM^i_{X/S})\stackrel{V_N}{\longrightarrow}H^j(X,\KK)
\end{equation}
are surjective.
Because of the relation $F_NV_N=N$ and the fact that multiplication by $N$ on
$H^j(X,\OM^i_{X/S})$ is injective, the map
$H^j(X,\OM^i_{X/S})\stackrel{V_N}{\longrightarrow}
H^j(X,\womn{N}^i_{X/S})$ is injective. Therefore the maps (\ref{eq:ker K}) are
in fact isomorphisms and the exact cohomology sequence of
(\ref{eq:left K}) breaks up into the short exact sequences (\ref{eq:Vi}).
\end{proof}

\

As a consequence of Propositions \ref{prop:wittsurjective2} and
\ref{prop:reduce ordinary} we have:

\begin{corollary}\label{cor:wittsurjective i}
Fix $N\geq 2$.
Assume in addition to the general hypotheses of this section that
for all prime numbers $p\leq N$ condition (\ref{eq:p}) holds, i.e.
$$
H^j(X,Z_p\OM_{X/S}^i)\rightarrow H^j(X,\OM_{X/S}^i)
\qquad\textrm{is an isomorphism for all}\quad i,j.
$$
Then for all $i,j$ and for every $n\leq N$ the sequence
\begin{equation}\label{eq:Vni}
0\rightarrow H^j(X,\OM^i_{X/S})\stackrel{V_n}{\longrightarrow}
H^j(X,\womn{n}^i_{X/S})\stackrel{\tr^n_{n-1}}{\longrightarrow}
H^j(X,\womn{n-1}^i_{X/S})\rightarrow 0
\end{equation}
is exact and, hence, the map, induced by truncation,
\begin{equation}\label{eq:surj Ni}
H^j(X,\womn{N}^i_{X/S})\rightarrow H^j(X,\OM^i_{X/S})
\end{equation}
is surjective.
\qed
\end{corollary}

\

Fix $N,i,j$.
Under the hypotheses of Proposition \ref{prop:wittsurjective1} (for $i=0$)
and Corollary \ref{cor:wittsurjective i} (for $i\geq 1$) the map, induced by
truncation,
$$
\tr^N_1: H^j(X,\womn{N}^i_{X/S})\rightarrow H^j(X,\OM^i_{X/S})
$$
is surjective. So we can pick a basis $\omega_1,\ldots,\omega_h$ for
$H^j(X,\OM^i_{X/S})$ and lift each $\omega_k$ to an
$\widetilde{\omega_k}$ in $H^j(X,\womn{N}^i_{X/S})$. For every $n\leq N$ we
denote the image of $\widetilde{\omega_k}$ under the standard truncation
$\tr^N_n:H^j(X,\womn{N}^i_{X/S})\rightarrow H^j(X,\womn{n}^i_{X/S})$
also by $\widetilde{\omega_k}$.

For every $n\leq N$ we then have an $h\times h$-matrix
$\mathrm{MAT}_{\widetilde{\omega}}(F_n)$ with entries in $A$ in which the
$k$-th column gives the coordinates of $F_n\widetilde{\omega_k}$ with respect
to the basis $\omega_1,\ldots,\omega_h$:
\begin{equation}\label{eq:mat F}
(F_n\widetilde{\omega_1},\ldots,F_n\widetilde{\omega_h})\,=\,
(\omega_1,\ldots,\omega_h)\cdot
\mathrm{MAT}_{\widetilde{\omega}}(F_n).
\end{equation}

\begin{lemma}\label{lemma:V generators}
Fix $N,i,j$.
Assume the hypotheses of Proposition \ref{prop:wittsurjective1} if $i=0$
and of Corollary \ref{cor:wittsurjective i} if $i\geq 1$. Let
$\omega_1,\ldots,\omega_h$ and
$\widetilde{\omega_1},\ldots,\widetilde{\omega_h}$ be as above.
Then for every $n\leq N$ and $\xi\in H^j(X,\womn{n}^i_{X/S})$ there are
unique elements $a_{l,k}\in A$ for $l=1,\ldots, n$ and $k=1,\ldots,h$ such that
\begin{equation}\label{eq:xi}
\xi=\sum_{k,l} V_l(\dul{a_{l,k}}\widetilde{\omega_k})
\end{equation}
where we have simplified the notation by writing $V_l$ instead of
$\tr^{nl+l-1}_n\circ V_l$.
\end{lemma}
\begin{proof}
Fix $n\leq N$ and $\xi\in H^j(X,\womn{n}^i_{X/S})$.
Since $\omega_1,\ldots,\omega_h$ is an $A$-basis of $H^j(X,\OM^i_{X/S})$
there are uniquely determined elements $a_{k,1}\in A$ such that:
$$
\tr^n_1 (\xi)=\sum_{k=1}^h a_{1,k}\omega_k.
$$
Then the truncation map
$\tr^n_2$ maps
$\xi-\sum_{k=1}^h \dul{a_{1,k}}\widetilde{\omega_k}$ into the kernel of
the map $\tr^2_1: H^j(X,\womn{2}^i_{X/S})\rightarrow H^j(X,\OM^i_{X/S})$.
{}From the exact sequence (\ref{eq:V0}) resp. (\ref{eq:Vni}) we see that there
are
$a_{2,1},\ldots,a_{2,h}\in A$ such that
$$
\tr^n_2(\xi-\sum_{k=1}^h \dul{a_{1,k}}\widetilde{\omega_k})=
V_2(\sum_{k=1}^h a_{2,k}\omega_k).
$$
Then the truncation map
$\tr^n_3$ maps
$\xi-\sum_{k=1}^h \dul{a_{1,k}}\widetilde{\omega_k}-
\sum_{k=1}^h V_2(\dul{a_{2,k}}\widetilde{\omega_k})$ into the kernel of
the map $\tr^3_2:H^j(X,\womn{3}^i_{X/S})\rightarrow H^j(X,\womn{2}^i_{X/S})$.
This kernel equals the image of
$V_3:H^j(X,\OM^i_{X/S})\rightarrow H^j(X,\womn{3}^i_{X/S})$. It is clear how
one can go on this way till one has reached (\ref{eq:xi}).
\end{proof}

\begin{lemma}\label{lemma:Fpnu surjective}
Fix $N,i,j$.
Assume the hypotheses of Proposition \ref{prop:wittsurjective1} if $i=0$
and of Corollary \ref{cor:wittsurjective i} if $i\geq 1$.
Let $p$ be a prime number, $p\leq N$. Assume that the map
\begin{equation}\label{eq:fp  sur}
F_p:F^*_pH^j(X,\womn{p}^i_{X/S})\rightarrow H^j(X,\OM^i_{X/S})
\qquad\textrm{is surjective.}
\end{equation}
Then for every $\nu\in\NN$ such that $p^\nu\leq N$
the map
\begin{equation}\label{eq:fpnu  sur}
F_{p^\nu}:F^*_{p^\nu}H^j(X,\womn{p^\nu}^i_{X/S})\rightarrow H^j(X,\OM^i_{X/S})
\qquad\textrm{is surjective.}
\end{equation}
\end{lemma}
\begin{proof}
Assume $\nu\geq 2$.
Let $\omega_1,\ldots,\omega_h$ and
$\widetilde{\omega_1},\ldots,\widetilde{\omega_h}$ be as above.
Recall that we also write $\widetilde{\omega_k}$ for the image of
$\widetilde{\omega_k}$ in $H^j(X,\womn{p^\nu}^i_{X/S})$ resp.
$H^j(X,\womn{p}^i_{X/S})$ under the truncation $\tr^N_{p^\nu}$ resp.
$\tr^N_{p}$. Then there are, according to the previous lemma,
unique elements $a_{l,k,m}\in A$ such that
$$
F_{p^{\nu-1}}\widetilde{\omega_m}=\sum_{k=1}^h\sum_{l=1}^p
V_l(\dul{a_{l,k,m,}}\widetilde{\omega_k}).
$$
Applying $F_p$ we see
$$
F_{p^{\nu}}\widetilde{\omega_m}\,=\,\sum_{k=1}^h\sum_{l=1}^p
F_pV_l(\dul{a_{l,k,m}}\widetilde{\omega_k})\\
\,=\,\sum_{k=1}^h (a_{1,k,m}^pF_p\widetilde{\omega_k}\:+\:pa_{p,k,m}\omega_k).
$$
This shows that with the notation from (\ref{eq:mat F}) we have
\begin{equation}\label{eq:fpnu}
\mathrm{MAT}_{\widetilde{\omega}}(F_{p^\nu})\equiv
\mathrm{MAT}_{\widetilde{\omega}}(F_p)\cdot
\mathrm{MAT}_{\widetilde{\omega}}(F_{p^{\nu-1}})^{(p)}
\;\bmod p
\end{equation}
where for a matrix $M=(a_{m,k})$ we denote $M^{(p^s)}=(a_{m,k}^{p^s})$.
Induction now shows
$$
\mathrm{MAT}_{\widetilde{\omega}}(F_{p^\nu})\equiv
\mathrm{MAT}_{\widetilde{\omega}}(F_p)\cdot
\mathrm{MAT}_{\widetilde{\omega}}(F_p)^{(p)}\cdot\ldots\cdot
\mathrm{MAT}_{\widetilde{\omega}}(F_p)^{(p^{\nu-1})}
\;\bmod p.
$$
Assumption (\ref{eq:fp  sur}) implies that
$\mathrm{MAT}_{\widetilde{\omega}}(F_p)$ is invertible $\bmod\, p$.
Thus $\mathrm{MAT}_{\widetilde{\omega}}(F_{p^\nu})$ is invertible $\bmod\, p$
for every $\nu$ such that $p^\nu\leq N$. This in turn implies
(\ref{eq:fpnu  sur}).
\end{proof}

\begin{lemma}\label{lemma:FN surjective}
Fix $N,i,j$.
Assume the hypotheses of Proposition \ref{prop:wittsurjective1} if $i=0$
and of Corollary \ref{cor:wittsurjective i} if $i\geq 1$.
Assume that for every prime number $p$ which divides $N$ the map
$$
F_p:F^*_pH^j(X,\womn{p}^i_{X/S})\rightarrow H^j(X,\OM^i_{X/S})
$$
is surjective.
Then the map
$$
F_N:F^*_NH^j(X,\womn{N}^i_{X/S})\rightarrow H^j(X,\OM^i_{X/S})
$$
is surjective.
\end{lemma}
\begin{proof}
Consider the prime decomposition $N=\prod_{l}p_l^{\nu_l}$.
Choose integers $c_l$ such that $\sum_l c_lNp_l^{-\nu_l}=1$.
Take $\xi\in H^j(X,\OM^i_{X/S})$.
{}From Lemma \ref{lemma:Fpnu surjective} we know that for every $l$ the map
$$
F_{p_l^{\nu_l}}:F^*_{p_l^{\nu_l}}H^j(X,\womn{{p_l^{\nu_l}}}^i_{X/S})
\rightarrow H^j(X,\OM^i_{X/S})
$$
is surjective. So there are
elements $a_{l,m}\in A$ and
$\alpha_{l,m}\in H^j(X,\womn{{p_l^{\nu_l}}}^i_{X/S})$
such that
$$
\xi= \sum_{m} a_{l,m} F_{p_l^{\nu_l}}\alpha_{l,m}.
$$
Then
$$
F_N\left(\sum_{l,m} c_l a_{l,m} V_{Np_l^{-\nu_l}}\alpha_{l,m}\right)
=
\sum_{l,m} c_lNp_l^{-\nu_l} a_{l,m} F_{p_l^{\nu_l}}\alpha_{l,m}
= \xi.
$$
So $F_N$ is surjective.
\end{proof}

\

\subsection{}

\

Theorem \ref{thm:U=conj} raised the issue of surjectivity of the map
$$
F_N:F^*_NH^j(X,\womn{N}^i_{X/S})\rightarrow H^j(X,\OM^i_{X/S}).
$$
Lemma \ref{lemma:FN surjective} reduced the problem to the surjectivity of
$$
F_p:F^*_pH^j(X,\womn{p}^i_{X/S})\rightarrow H^j(X,\OM^i_{X/S})
$$
for prime numbers $p$ dividing $N$. Recall that in the circumstances of
Proposition \ref{prop:wittsurjective1} and  Corollary \ref{cor:wittsurjective
i} one has short exact sequences
$$
0\rightarrow H^j(X,\OM^i_{X/S})\stackrel{V_n}{\longrightarrow}
H^j(X,\womn{n}^i_{X/S})\stackrel{\tr^n_{n-1}}{\longrightarrow}
H^j(X,\womn{n-1}^i_{X/S})\rightarrow 0
$$
{}From these sequences and the relations
$$
F_pV_p=p\,,\qquad F_pV_n=V_nF_p\quad\textrm{if}\;n<p
$$
one deduces that there is a map
\begin{equation}\label{eq:F mod p}
\oF_p:H^j(X,\OM^i_{X/S})\q{p}\rightarrow H^j(X,\OM^i_{X/S})\q{p}
\end{equation}
making the following square commutative
$$
\begin{array}{ccc}
H^j(X,\womn{p}^i_{X/S})&\stackrel{F_p}{\longrightarrow} &H^j(X,\OM^i_{X/S})\\
[1ex]
(\bmod p)\circ\tr^p_1\downarrow&&\downarrow \bmod p\\[1ex]
H^j(X,\OM^i_{X/S})\q{p}&\stackrel{\oF_p}{\longrightarrow}&
H^j(X,\OM^i_{X/S})\q{p}
\end{array}
$$
Set
$$
\oF_p^*H^j(X,\OM^i_{X/S})=(A/pA)\otimes_A H^j(X,\OM^i_{X/S})
$$
where $A/pA$ is an $A$-module via the ring homomorphism $A\rightarrow A/pA,\,
{a\mapsto a^p\bmod p}$.
Then $\oF_p$ induces the $A$-linear map
$$
\oF_p:\oF_p^*H^j(X,\OM^i_{X/S})\rightarrow H^j(X,\OM^i_{X/S})\,,\qquad
\oF_p (a\otimes\alpha):=a\cdot\oF_p\alpha.
$$
Now it is clear that
$$
F_p:F^*_pH^j(X,\womn{p}^i_{X/S})\rightarrow H^j(X,\OM^i_{X/S})\qquad\textrm{is
surjective}
$$
if and only if
$$
\oF_p:\oF_p^*H^j(X,\OM^i_{X/S})\rightarrow H^j(X,\OM^i_{X/S})\qquad\textrm{is
surjective}.
$$
To prove surjectivity for the latter $\oF_p$ we need all hypotheses from
Setting \ref{setting}.

\

\textbf{So from now on }
$f:X\rightarrow S=\Spec A$ is a smooth projective morphism of relative
dimension $\rd$, such that the ring $A$ is \'etale over a polynomial ring
$\ZZ\half[x_1,\ldots,x_r]$ and  $H_{DR}^m(X/S)$ and
$H^{m-i}(X,\OM_{X/S}^i)$ are free $A$-modules for all $m,i$.

Let $\OM^\rd_{X/S,\mathrm{log}}$ be the sheaf of \emph{logarithmic
$\rd$-forms}. It is locally generated by sections of the form
\begin{equation}\label{eq:log forms}
\frac{d u_1}{u_1}\wedge\ldots\wedge \frac{d u_\rd}{u_\rd}.
\end{equation}
We assume that
$H^\rd (X,\OM^\rd _{X/S})$ is a free $A$-module of rank $1$
and that we can choose $\varpi\in H^\rd (X,\OM^\rd _{X/S,\log})$ so that its
image $\varpi_1$ in $H^\rd (X,\OM^\rd _{X/S})$ is an $A$-basis of the latter.
Using the product map
$$
H^j(X,\OM^i_{X/S})\times H^{\rd -j}(X,\OM^{\rd -i}_{X/S})
\rightarrow H^\rd (X,\OM^\rd _{X/S})
$$
in the Hodge cohomology algebra
and the chosen $\varpi$ we define the bilinear pairing
\begin{equation}\label{eq:bilin}
\langle\:,\:\rangle:\;
H^j(X,\OM^i_{X/S})\times H^{\rd -j}(X,\OM^{\rd -i}_{X/S})
\rightarrow A\qquad\textrm{s.t.}\qquad
\alpha\cdot\beta=\langle\alpha,\beta\rangle\varpi_1.
\end{equation}
We assume that this pairing is non-degenerate for all $i,j$.

\begin{theorem}\label{thm:Hodge-Tate}
Fix $N\geq 2$.
Assume in addition to the above hypotheses
that $H^j(X,Z_p\OM_{X/S}^i)\rightarrow H^j(X,\OM_{X/S}^i)$
is an isomorphism for all $i,j$ and for all prime numbers $p\leq N$.
Then for all $i,j$ and for every prime number $p\leq N$ the map
$$
\oF_p:\oF_p^*H^j(X,\OM^i_{X/S})\rightarrow H^j(X,\OM^i_{X/S})\q{p}
$$
is an isomorphism.
\end{theorem}
\begin{proof}
A section of the form (\ref{eq:log forms}) can be lifted to a section
$$
(\dul{u_1}\cdot\ldots\cdot \dul{u_\rd})^{-1}
d \dul{u_1}\wedge\ldots\wedge
d \dul{u_\rd}
$$
of $\womn{n}^\rd_{X/S}$ for every $n$.
Thus we get a morphism of sheaves
$$
\OM^\rd_{X/S,\mathrm{log}}\rightarrow\womn{n}^\rd_{X/S}
$$
which composes with the truncation $\tr^n_1: \womn{n}^\rd_{X/S}\rightarrow
\OM^\rd_{X/S}$ to the inclusion $\OM^\rd_{X/S,\mathrm{log}}\subset
\OM^\rd_{X/S}$.
Let $\varpi_n\in H^\rd(X,\womn{n}^\rd_{X/S})$ denote the image of $\varpi$
under the map
$H^\rd(\OM^\rd_{X/S,\mathrm{log}})\rightarrow H^\rd(\womn{n}^\rd_{X/S})$.
Using the formulas (\ref{eq:FVrel}), first on sections and then on cohomology,
one checks
\begin{equation}\label{eq:varpi fixed}
F_l\varpi_{nl}=\varpi_n\qquad\textrm{for all}\quad n,l.
\end{equation}
Take a basis $\omega_1,\ldots,\omega_h$ for
$H^j(X,\OM^i_{X/S})$ and a basis $\eta_1,\ldots,\eta_h$ for
$H^{\rd-j}(X,\OM^{\rd-i}_{X/S})$ such that the matrix
$$
\left(\langle\omega_k,\eta_m\rangle\right)_{m,k}
$$
is invertible over $A$; for this we use the assumption that the pairing
$\langle\,,\,\rangle$ is non-degenerate.
Lift each $\omega_k$ to an
$\widetilde{\omega_k}$ in $H^j(X,\womn{N}^i_{X/S})$
and lift each $\eta_m$ to an
$\widetilde{\eta_m}$ in $H^{\rd-j}(X,\womn{N}^{\rd-i}_{X/S})$.
Then, as in Lemma \ref{lemma:V generators}, there are unique elements
$a_{l,k,m}\in A$ such that in $H^{\rd}(X,\OM^{\rd}_{X/S})$
\begin{equation}\label{eq:IP}
\widetilde{\omega_k}\cdot\widetilde{\eta_m}=
\sum_lV_l(\dul{a_{l,k,m}}\varpi_N).
\end{equation}
Note $a_{1,k,m}=\langle\omega_k,\eta_m\rangle$.
Let $p$ be a prime number, $p\leq N$. Apply the operator $F_p\circ\tr^N_p$
to (\ref{eq:IP}):
$$
(F_p\circ\tr^N_p\widetilde{\omega_k})\cdot
(F_p\circ\tr^N_p\widetilde{\eta_m})=
(a_{1,k,m}^p+pa_{p,k,m})\varpi_1.
$$
Taken modulo $p$ this relation reads
$$
\langle\oF_p(\omega_k\bmod p),\oF_p(\eta_m\bmod p)\rangle
=
\langle\omega_k\bmod p\,,\,\eta_m\bmod p\rangle^p.
$$
In matrix form this is
$$
\mathrm{MAT}_{{\omega}}(\oF_{p})^{\mathrm{transpose}}
\cdot
\left(\langle\omega_k,\eta_m\rangle\right)
\cdot
\mathrm{MAT}_{{\eta}}(\oF_{p})
\equiv
\left(\langle\omega_k,\eta_m\rangle^p\right)\;\bmod p.
$$
So the matrix $\mathrm{MAT}_{{\omega}}(\oF_{p})$
is invertible over $A/pA$.
\end{proof}

\begin{remark}
For $(i,j)\neq (\frac{\rd}{2},\frac{\rd}{2})$ one can, after an initial choice
of bases $\{\omega_1,\ldots,\omega_h\}$ for
$H^j(X,\OM^i_{X/S})$ and $\{\eta_1,\ldots,\eta_h\}$ for
$H^{\rd-j}(X,\OM^{\rd-i}_{X/S})$, change the latter by means of the matrix
$\left(\langle\omega_k,\eta_m\rangle\right)^{-1}$ and thus obtain a basis
$\{\zeta_1,\ldots,\zeta_h\}$ for $H^{\rd-j}(X,\OM^{\rd-i}_{X/S})$ such that
$$
\langle\omega_k,\zeta_m\rangle=\delta_{k,m}\qquad\textrm{(Kronecker's delta)}
$$
Then the above proof shows
\begin{equation}\label{eq:Hasse-Witt inverse}
\mathrm{MAT}_{{\zeta}}(\oF_{p})
=
\mathrm{MAT}_{{\omega}}(\oF_{p})^{\mathrm{transpose\; inverse}}
\end{equation}
\end{remark}

\

\section*{APPENDIX: Formal groups in the background}
Behind the structures in Proposition \ref{prop:wittsurjective1} and
Corollary \ref{cor:wittsurjective i} one can see formal groups.
In this appendix we discuss the formal groups for Proposition
\ref{prop:wittsurjective1}, leaving the ones for Corollary
\ref{cor:wittsurjective i} for future investigations.
The formal groups for Proposition \ref{prop:wittsurjective1} are well known
and were introduced by
Artin and Mazur \cite{AM} as a generalization of the classical
construction/definition of the \emph{formal Picard group}: the latter is based
on deformation of $H^1(X,\cO_X^*)$; the new ones are  based on
deformation of $H^j(X,\cO_X^*)$ for any $j$. Despite the abstract appearance of
their definition these formal groups are often
very manageable objects which can be used for concretely computing the action
of Frobenius operators on the cohomology of sheaves of generalized Witt
vectors.

We place our discussion of the Artin-Mazur formal groups in the geometric
setting of this paper, i.e. $\XX$, $\bS$ and $\bA$ are as in Setting
\ref{setting}. Unlike ``ordinariness'' Artin-Mazur formal groups require no
further restrictive conditions; everything works with the original $\XX$, $\bS$
and $\bA$ .

The Artin-Mazur formal group $H^j (\XX,\hat{\GG}_{m,\XX})$ is defined as
follows ($m$ in the notation
$\hat{\GG}_{m,\XX}$ refers to the \emph{multiplicative} formal group).
Let $\nar$ denote the category of
nil-$\bA$-algebras, i.e. associative commutative $\bA$-algebras without unit
element in which every element is nilpotent.
For a nil-$\bA$-algebra $\nil$ we have the sheaf of nilalgebras
$\cO_\XX\otimes_\bA \nil$ on $\XX$. For local sections $x$ and $y$ of
$\cO_\XX\otimes_\bA \nil$ we set
$$
x\star y=x+y-xy.
$$
The sheaf $\cO_\XX\otimes_\bA \nil$ with the binary operation $\star$ on its
local sections is a sheaf of abelian groups on $\XX$, which we denote by
$\hat{\GG}_{m,\XX}(\nil)$.
One defines $H^j(\XX,\hat{\GG}_{m,\XX})$ to be the functor
\begin{eqnarray}
\nonumber
H^j(\XX,\hat{\GG}_{m,\XX})\::\;\nar&\rightarrow&
\mathfrak{Abelian\; groups}
\\
\label{eq:AM}
&&
\\
\nonumber
\nil&\mapsto&H^j(\XX,\hat{\GG}_{m,\XX}(\nil)).
\end{eqnarray}
In the category $\nar$ one has in particular the nilalgebras
$t\bA[t]/(t^{n+1})$ for $n\geq 1$. Noticing that
$$
\hat{\GG}_{m,\XX}(t\bA[t]/(t^{n+1}))=1+t\OO_\XX[t]/(t^{n+1})=\won{n}_\XX
$$
we recover the Witt vector cohomology groups
$$
H^j(X,\hat{\GG}_{m,\XX}(t\bA[t]/(t^{n+1})))=H^j(X,\won{n}_\XX).
$$
According to Proposition \ref{prop:wittsurjective1} the truncation maps
$$
H^j(X,\won{n}_\XX)\rightarrow H^j(X,\won{n-1}_\XX)
$$
are surjective. So the limit $\lim_{\leftarrow n}H^j(X,\won{n}_\XX)$
contains elements which project onto a basis of $H^j(X,\OO_\XX)$
as $\bA$-module.
A convenient notation for an element of
$\lim_{\leftarrow n}H^\rd(X,\won{n}_\XX)$ is $\gamma(t)$.
For
$$
\gamma(t)\in \lim_{\leftarrow n}H^\rd(X,\won{n}_\XX)
$$
we construct, as follows, a functorial homomorphism
$$
\gamma: \nil\rightarrow H^j(\XX,\hat{\GG}_{m,\XX}(\nil))
$$
from on the one side the functor $\nar\rightarrow\mathfrak{Sets}$ which assigns
to a nilalgebra $\nil$ the underlying set $\nil$ to on the other side the
functor
which assigns to $\nil$ the set underlying the group
$H^j(\XX,\hat{\GG}_{m,\XX}(\nil))$.
For the construction of $\gamma$ take a nilalgebra $\nil$ and an element
$u\in\nil$. There is then an $n\in\NN$ and an algebra homomorphism
$$
g:t\bA[t]/(t^{n+1})\rightarrow\nil\,,\qquad g(t)=u.
$$
This $g$ induces a group homomorphism
$$
g_*: H^j(X,\won{n}_\XX)\rightarrow H^j(\XX,\hat{\GG}_{m,\XX}(\nil)).
$$
We define
$$
\gamma(u)=g_*(\gamma(t))\in H^\rd(\XX,\hat{\GG}_{m,\XX}(\nil)).
$$
and obtain thus a map of sets
$$
\gamma: \nil\rightarrow H^j(\XX,\hat{\GG}_{m,\XX}(\nil))\,,\quad u\mapsto
\gamma(u).
$$
The functoriality of this construction is obvious.

Now fix elements $\gamma_{j,1}(t),\ldots,\gamma_{j,h^{0j}}(t)$ in
$\lim_{\leftarrow n}H^j(X,\won{n}_\XX)$
such that their images in $H^j(X,\OO_\XX)$ constitute an $\bA$-module basis.
Then the above construction gives a functorial map
\begin{equation}\label{eq:formal coord}
\begin{array}{cccc}
(\gamma_{j,1},\ldots,\gamma_{j,h^{0j}}):& \nil^{\times h^{0j}}&\rightarrow&
H^j(\XX,\hat{\GG}_{m,\XX}(\nil))\\
&(u_1,\ldots,u_{h^{0j}})&\mapsto& \sum_{k=1}^{h^{0j}}\gamma_k(u_k)\end{array}
\end{equation}

\begin{proposition}\label{prop:formal coord}
(\ref{eq:formal coord}) is a functorial bijection for every $j$.
\end{proposition}
\begin{proof}
First note that by a simple inductive limit argument it suffices to show
bijectivity only for nilalgebras which are finitely generated over $\bA$. For
the finitely generated nilalgebras
one proceeds by induction along so-called \emph{small extensions}: a small
extension is a surjective algebra homomorphism $\nil\twoheadrightarrow\nil'$
with kernel generated by a single element $\epsilon$ such that
$\epsilon\nil=0$. For every finitely generated nilalgebra $\nil$ the trivial
map $\nil\rightarrow 0$ is the composite of finitely many small extensions.

For the induction step ones notes that
$$
\hat{\GG}_{m,\XX}(\epsilon\bA))=1+\epsilon\bA\otimes\OO_\XX
\stackrel{\log}{\simeq}\epsilon\bA\otimes\OO_\XX
$$
and hence
\begin{equation}\label{eq:tangent space}
H^j(\XX,\hat{\GG}_{m,\XX}(\epsilon\bA))\stackrel{\log}{\simeq}
H^j(\XX,\OO_\XX\otimes\epsilon\bA)=
H^j(\XX,\OO_\XX)\otimes\epsilon\bA\stackrel{\mathrm{basis}}{\simeq}
(\epsilon\bA)^{\times h^{0j}}.
\end{equation}
For a small extension $\epsilon\bA\hookrightarrow\nil\twoheadrightarrow\nil'$
the sequence of sheaves of groups on $\XX$
\begin{equation}\label{eq:exact sheaves}
0\rightarrow\hat{\GG}_{m,\XX}(\epsilon\bA)
\rightarrow\hat{\GG}_{m,\XX}(\nil))
\rightarrow\hat{\GG}_{m,\XX}(\nil')\rightarrow 0
\end{equation}
is obviously exact.
There is, by functoriality, a commutative diagram
\begin{equation}\label{eq:coord ladder}
\begin{array}{ccccc}
(\epsilon\bA)^{\times h^{0j}}
&\hookrightarrow&\nil^{\times h^{0j}}&\twoheadrightarrow&
(\nil')^{\times h^{0j}}\\
\downarrow&&\downarrow&&\downarrow\\
H^j(\XX,\hat{\GG}_{m,\XX}(\epsilon\bA))&\rightarrow&
H^j(\XX,\hat{\GG}_{m,\XX}(\nil))&\rightarrow&
H^j(\XX,\hat{\GG}_{m,\XX}(\nil'))
\end{array}
\end{equation}
in which the vertical maps are given by (\ref{eq:formal coord})
and the bottom row is part of the exact sequence of cohomology groups
associated with (\ref{eq:exact sheaves}):
$$
\begin{array}{rlcll}
0\rightarrow&
H^0(\XX,\hat{\GG}_{m,\XX}(\epsilon\bA))\rightarrow
& H^0(\XX,\hat{\GG}_{m,\XX}(\nil))&\rightarrow
H^0(\XX,\hat{\GG}_{m,\XX}(\nil'))&\rightarrow\\
\rightarrow&
H^1(\XX,\hat{\GG}_{m,\XX}(\epsilon\bA))\rightarrow
&\ldots\quad\ldots\quad\ldots&
\rightarrow
H^{\rd -1}(\XX,\hat{\GG}_{m,\XX}(\nil'))
&\rightarrow\\
\rightarrow&
H^\rd(\XX,\hat{\GG}_{m,\XX}(\epsilon\bA))\rightarrow
&H^\rd(\XX,\hat{\GG}_{m,\XX}(\nil))&\rightarrow
H^\rd(\XX,\hat{\GG}_{m,\XX}(\nil'))
&\rightarrow 0
\end{array}
$$
The proposition can be proved by a double induction with inside an induction
for $j=0,\ldots,\rd$ an induction along small extensions.
The result for $j-1$ provides surjectivity of
$
H^{j-1}(\XX,\hat{\GG}_{m,\XX}(\nil))\rightarrow
H^{j-1}(\XX,\hat{\GG}_{m,\XX}(\nil'))
$
and thus injectivity of
$H^j(\XX,\hat{\GG}_{m,\XX}(\epsilon\bA))\rightarrow
H^j(\XX,\hat{\GG}_{m,\XX}(\nil))$. Knowing this injectivity and assuming that
in (\ref{eq:coord ladder}) the left and right vertical maps are bijective one
can conclude that the middle map is bijective too.
\end{proof}

\

A functorial bijection such as (\ref{eq:formal coord}) is called a
\emph{coordinatization of the functor $H^j(\XX,\hat{\GG}_{m,\XX})$.}
A functor which allows such a coordinatization (i.e. a functorial bijection to
the set valued functor $\nil\mapsto \nil^{\times h}$ for some $h$)
\emph{is called a formal group}. So in the circumstances of Setting
\ref{setting} the Artin-Mazur functors are indeed formal groups.

There is a converse to the above result in that every coordinatization of the
functor $H^j(\XX,\hat{\GG}_{m,\XX})$ provides a set of elements in
$\lim_{\leftarrow n}H^j(X,\won{n}_\XX)$
such that their images in $H^j(X,\OO_\XX)$ constitute an $\bA$-module basis.
To see this one just has to evaluate the coordinatization at the
projective system of nilalgebras
$\left\{t\bA[t]/(t^{n+1})\right\}_{n\geq 1}$ and to observe that functoriality
implies that the coordinatization
$$
(t\bA[t]/(t^2))^{h^{0j}}\longrightarrow H^j(X,\OO_X)
$$
is not only a bijective map of sets but even an $\bA$-linear isomorphism.

\

We are mainly interested in applying the theory to families of Calabi-Yau
varieties. Since we do not want to blur the picture with unnecessary
technicalities from formal group theory we assume henceforth
\begin{equation}\label{eq:dim=1}
H^\rd(\XX,\OO_\XX)\simeq\bA
\end{equation}
and we look  only at the
Artin-Mazur formal group $H^\rd(\XX,\hat{\GG}_{m,\XX})$.
A coordinatization of this $1$-dimensional formal group, i.e. a functorial
bijection
$$
\gamma: \nil\stackrel{1:1}{\longrightarrow} H^\rd(\XX,\hat{\GG}_{m,\XX}(\nil)),
$$
provides a \emph{formal group law}, i.e. a two variable power series
$G(t_1,t_2)\in \bA[[t_1,t_2]]$: if $y_1,y_2$ are elements in some nilalgebra
$\nil$, then
$$
\gamma(y_1)+\gamma(y_2)= \gamma(G(y_1,y_2))\qquad\textrm{in the group}\;
H^\rd(\XX,\hat{\GG}_{m,\XX}(\nil)).
$$
Since the ring $\bA$ is flat over $\ZZ$ (i.e. the canonical map
$\bA\rightarrow \bA\otimes \QQ$ is injective), there is a power series
$$
\ell(t)\:=\:\sum_{m\geq 1}
\frac{a_m}{m}\,t^m\;\in (\bA\otimes\QQ)[[t]]
$$
in one variable $t$ with all $a_m\in \bA$ and $a_1=1$, such that
$$
G(t_1,t_2)=\ell^{-1}(\ell(t_1)+\ell(t_2)).
$$
The fact that the coefficients of the power series $G(t_1,t_2)$ have no
denominators, puts a very strong structure on the sequence $\{a_m\}_{m\geq 1}$.

\

In \cite{S3} Theorem 1 it is shown how a concrete logarithm for the Artin-Mazur
formal group $H^\rd(\XX,\hat{\GG}_{m,\XX})$ can be obtained from the defining
equations when $\XX$ is a complete intersection  in $\PP_\bA^N$. Proposition
\ref{logarithm} gives the result for Calabi-Yau varieties of dimension $\rd$
(conditions $d_1+\ldots+d_r\,=\,N+1$ and $N-s-1=\rd$).

\begin{proposition}\label{logarithm}
Let $\bA$ be as in Setting \ref{setting}. Let $P_1,\ldots,P_s$ be a
\emph{regular sequence} of homogeneous polynomials in $\bA[Z_0,\ldots,Z_N]$ of
degrees $d_1,\ldots,d_s$ with  $d_1+\ldots+d_s\,=\,N+1$ and $N-s=\rd$.
Let $\XX$ be the subscheme of $\PP_\bA^{N}$ defined by the ideal
$(P_1,\ldots,P_s)$. Assume $\XX$ is flat over $\bA$.
Then the power series
$$
\ell (t)\:=\:\sum_{m\geq 1}\frac{a_m}{m}\,t^m
$$
with
$$
a_m\::=\:\textrm{\textup{ coefficient of }} (Z_0\cdot\ldots\cdot Z_N)^{m-1}
\hspace{3mm}\mathrm{ in }\hspace{3mm}
(P_1\cdot\ldots\cdot P_s)^{m-1}
$$
is a logarithm for the Artin-Mazur formal group $H^\rd(\XX,\hat{\GG}_{m,\XX})$.
\qed
\end{proposition}

\

Well-known examples in which the conditions of this Proposition
\ref{logarithm} are satisfied, are:
$\diamond$
quintic hypersurfaces in $\PP^4$
$\diamond$
intersections of two cubic hypersurfaces in $\PP^5$
$\diamond$
intersections of four quadrics in $\PP^7$.

\cite{S3} Theorem 2 gives a similar result for branched double coverings of
$\PP_\bA^{N}$ and applies for instance to double coverings of $\PP_\bA^3$
branched along a hypersurface of degree $8$.

Because of the enormous number of variables it is not practical to write out
the $a_m$'s explicitly for the general families in these examples. There are,
however, manageable subfamilies, like:

\begin{example}\label{ex:quintic}
For the family of Calabi-Yau hypersurfaces in $\PP^4$
with equation
$$
Z_0Z_1Z_2Z_3Z_4 -x (Z_0^5+Z_1^5+Z_2^5+Z_3^5+Z_4^5)=0
$$
one has
$$
a_m=\sum_{j\geq 0} \; (-1)^j\frac{(5j)!}{(j!)^5}
\left(\begin{array}{c} m-1\\5j \end{array}\right)\:x^{5j}.
$$
Note that $a_m$ is a polynomial of degree $\leq m-1$.
\end{example}

\

The group
$$
\lim_{\leftarrow n}H^\rd(X,\hat{\GG}_{m,\XX}(t\bA[t]/(t^{n+1})))=
\lim_{\leftarrow n}H^\rd(X,\won{n}_\XX)
$$
is called the \emph{module of curves on} or the \emph{Cartier-Dieudonn\'e
module} of the formal group $H^\rd(X,\hat{\GG}_{m,\XX})$. It comes equipped
with operators $\dul{a}$, $V_k$ and $F_k$ for $a\in\bA$, $k\in \NN$, defined
as acting on elements
$\gamma (t)\in\lim_{\leftarrow n}H^\rd(X,\hat{\GG}_{m,\XX}(t\bA[t]/(t^{n+1})))$
by
\begin{equation}\label{eq:operators on curves}
\begin{array}{rcll}
\dul{a}\gamma (t)&=&\gamma(at)
&\textrm{i.e. induced by substution}\quad t\mapsto at\\
V_k\gamma (t)&=&\gamma(t^k)
&\textrm{i.e. induced by substution}\quad t\mapsto t^k\\
F_k\gamma (t)&=&\boxplus_{l=0}^{k-1}\gamma(\zeta^l t^{1/k})
\quad& \textrm{with }\;\zeta\;\textrm{ primitive }\;k\textrm{-th root of unity}
\end{array}
\end{equation}
where $\boxplus$ refers to the addition in the formal group, which corresponds
with the standard addition in
$\lim_{\leftarrow n}H^\rd(X,\won{n}_\XX)$.
The operators $F_k$ and $V_k$ correspond with the earlier defined
Frobenius and Verschiebung operators on
$\lim_{\leftarrow n}H^\rd(X,\won{n}_\XX)$; the operator $\dul{a}$ corresponds
with multiplication by the element
$\dul{a}\in\underline{{\mathcal {W}}\bA}$.

If $\gamma(t)$ projects onto an $\bA$-basis of $H^\rd(\XX,\OO_\XX)$ it gives
rise to a coordinatization, a formal group law and a logarithm
$$
\ell(t)\:=\:\sum_{m\geq 1}
\frac{a_m}{m}\,t^m
$$
for this formal group law. Since the coordinatization maps $t$ to $\gamma(t)$
we can now concretely compute
\begin{eqnarray*}
F_k\gamma (t)&=&\boxplus_{l=0}^{k-1}\gamma(\zeta^l t^{1/k})\;
=\;\ell^{-1}\left(\sum_{l=0}^{k-1}\ell(\zeta^l t^{1/k})\right)\\
&=&\ell^{-1}\left(\sum_{l=0}^{k-1}\sum_{m\geq 1}
\frac{a_m}{m}(\zeta^l t^{1/k})^m\right)\\
&=&\ell^{-1}\left(a_kt+\textrm{ higher order terms}\right)\\
&=&a_kt\;\bmod t^2.
\end{eqnarray*}
Looking back at (\ref{eq:mat F}) and the discussion  preceding it we see:
we have taken a basis $\overline{\gamma}$ for $H^\rd(\XX,\OO_\XX)$,
lifted it to an element
$\gamma (t)\in\lim_{\leftarrow n}
H^\rd(X,\hat{\GG}_{m,\XX}(t\bA[t]/(t^{n+1})))=
\lim_{\leftarrow n}H^\rd(X,\won{n}_\XX)
$ and computed the image of $F_k\gamma (t)$ in $H^\rd(\XX,\OO_\XX)$.
The conclusion is
\begin{equation}\label{eq:frob concrete}
\mathrm{MAT}_{\gamma (t)}(F_k)=a_k.
\end{equation}
\emph{
This makes a direct link between the Frobenius operators
$$
F_k:H^\rd(X,\won{k}_\XX)\rightarrow H^\rd(X,\OO_\XX)
$$
and the coefficients
of the logarithm for a formal group law of the Artin-Mazur formal group
$H^\rd(X,\hat{\GG}_{m,\XX})$.}

\

Recall the following special case of (\ref{eq: G-M nilpotent}):
$$
\nabla\wtF_kH^\rd(\XX,\won{k}_\XX)\,=\,0.
$$
This means here concretely
$$
\nabla(a_k\overline{\gamma})\equiv 0 \,\bmod k.
$$
Let $\omega\in H^0(\XX,\OM^\rd_{\XX/\bS})$ be such that
$$
\langle\omega,\overline{\gamma}\rangle=1.
$$
The Gauss-Manin connection $\nabla$ on $H^\rd_{DR}(\XX/\bS)$ induces an action
of the ring of differential operators $\mathrm{Diff}(\bS/\ZZ)$
(or rather of the subring generated by the derivations).
{}From the above discussion we derive the following result
which is a special case of \cite{S2} thm. 4.6:

\begin{proposition}
With the above hypotheses and notations one has:\\
if the differential operator $L\in\mathrm{Diff}(\bS/\ZZ)$ is such that
$$
L\omega=0,
$$
then
$$
La_k\equiv 0\bmod k\bA.
$$
\qed
\end{proposition}

\begin{example}\textup{(continuation of Example \ref{ex:quintic})}
For an appropriate choice of the nowhere vanishing holomorphic $3$-form
$\omega$ the \emph{Picard-Fuchs operator}, i.e. the differential operator
annihilating $\omega$, is
$$
L=\theta^4-5^5x^5(\theta+1)(\theta+2)(\theta+3)(\theta+4)\qquad\quad
\textrm{with} \quad \theta=x\frac{d}{dx}.
$$
One easily checks that for
\begin{equation}\label{eq:ak}
a_k=\sum_{j\geq 0} \; (-1)^j\frac{(5j)!}{(j!)^5}
\left(\begin{array}{c} k-1\\5j \end{array}\right)\:x^{5j}
\end{equation}
as in Example \ref{ex:quintic}
\begin{equation}\label{eq:ak diff}
La_k\equiv 0\bmod k\ZZ[x]
\end{equation}
On the other hand one has the power series
\begin{equation}\label{eq:f}
f(x)=\sum_{j\geq 0} \; \frac{(5j)!}{(j!)^5}\:x^{5j}
\end{equation}
which satisfies
\begin{equation}\label{eq:f diff}
Lf=0.
\end{equation}
\end{example}

\

The striking resemblance between (\ref{eq:ak}) and (\ref{eq:f})
and between (\ref{eq:ak diff}) and (\ref{eq:f diff}) holds in great generality
and shows where and how the ``holomorphic solution of the Picard-Fuchs equation
near the large complex structure limit''
shows up in our ``ordinary limit''.

\


\begin{thebibliography}{99}
%
\bibitem{AM}
Artin, M., and B. Mazur, \textit{Formal groups arising from algebraic
varieties}, Ann. scient. \'Ec. Norm. Sup. 4e s\'erie, \textbf{10} (1977),
87--132.
%
\bibitem{B}
Bourbaki, N. \textit{\'El\'ements de math\'ematique: Alg\`ebre commutative,
chapitres 8 et 9}, Masson, Paris, 1983.
%
\bibitem{C}
Cartier, P. {\em Groupes formels associ\'{e}s aux anneaux de Witt
g\'{e}n\'{e}ralis\'{e}s. } resp. {\em Modules associ\'{e}s \`{a} un
groupe formel commutatif. Courbes typiques. } C.R. Acad. Sc. Paris 265
(1967) 49-52 resp. 129-132
%
\bibitem{D}
Deligne, P., \textit{Local behavior of Hodge structures at infinity}, Mirror
Symmetry II, B. Greene and S-T. Yau (eds.), AMS/IP Studies in Advanced Math.
vol. 1, Amer. Math. Soc. and International Press, 1997, pp. 683--699.
%
\bibitem{DI}
Deligne, P., L. Illusie, \textit{Rel\`evements modulo $p^2$ et d\'ecomposition
du complexe de De Rham}, Invent. Math. \textbf{89} (1987) 247--270.
%
\bibitem{H}
Hazewinkel, M. {\em Formal groups and applications. }
New York: Academic Press 1978
%
\bibitem{I}
Illusie, L., \textit{Complexe de De Rham-Witt et cohomologie cristalline},
Ann. scient. \'{E}c. Norm. Sup. 4e s\'erie, \textbf{12} (1979) 501--661.
%
\bibitem{I2}
Illusie, L., \textit{Ordinarit\'e des intersections compl\`etes g\'en\'erales},
The Grothendieck Festschrift volume II, P. Cartier et al. (eds.), Progress in
Math. vol. 87, Birkh\"auser, Boston, 1990, pp. 375--405.
 %
\bibitem{IR}
Illusie, L., and M. Raynaud, \textit{Les suites spectrales associ\'ees au
complexe de De Rham-Witt}, Publ. Math. IHES \textbf{57} (1983) 73--212.
%
\bibitem{K}
Katz, N. \textit{Nilpotent connections and the monodromy theorem: applications
of a result of Turrittin}, Publ. Math. IHES \textbf{39} (1970) 175--232.
%
\bibitem{M}
Morrison, D., \textit{Mirror symmetry and rational curves on quintic
threefolds: a guide for mathematicians}, J. of the Amer. Math. Soc. \textbf{6}
(1993), 223--247.
%
\bibitem{S1}
Stienstra, J., \textit{Ordinary Calabi-Yau-3 crystals}, these proceedings.
%
\bibitem{S2}
Stienstra, J., \textit{The generalized De Rham-Witt complex and congruence
differential equations}, Arithmetic Algebraic Geometry, G. van der Geer, F.
Oort, J. Steenbrink (eds.), Progress in Math. vol. 89, Birkh\"auser, Boston,
1991, pp. 337--358.
%
\bibitem{S3}
Stienstra, J., \textit{Formal group laws arising from algebraic varieties},
Amer. J. Math. \textbf{109} (1987) 907--925

%
\end{thebibliography}
\end{document}